\documentclass{elsarticle}

\usepackage[colorlinks=true, urlcolor=black, citecolor=black, linkcolor=black, hyperfootnotes=true]{hyperref}
\usepackage{amssymb}
\usepackage{amsmath}
\usepackage{aliascnt}

\numberwithin{equation}{section}

\newtheorem{thm}{Theorem}[section]

\newaliascnt{prp}{thm}
\newtheorem{prp}[prp]{Proposition}
\aliascntresetthe{prp}

\newaliascnt{cor}{thm}
\newtheorem{cor}[cor]{Corollary}
\aliascntresetthe{cor}

\newaliascnt{dfn}{thm}
\newdefinition{dfn}[dfn]{Definition}
\aliascntresetthe{dfn}

\newaliascnt{xpl}{thm}
\newdefinition{xpl}[xpl]{Example}
\aliascntresetthe{xpl}

\newaliascnt{rmk}{thm}
\newdefinition{rmk}[rmk]{Remark}
\aliascntresetthe{rmk}

\newdefinition{ass}{Assumption}

\newproof{proof}{Proof}

\author{Tristan Bice\tnoteref{grantinfo}}
\tnotetext[grantinfo]{The author is supported by the GA\v{C}R project EXPRO 20-31529X and RVO: 67985840.}
\address{Institute of Mathematics\\
Czech Academy of Sciences\\
\v{Z}itn\'a 25, 115 67 Prague, Czech Republic}
\ead{Tristan.Bice@gmail.com}

\title{An Algebraic Approach to the Weyl Groupoid}

\begin{document}

\begin{abstract}
We unify the Kumjian-Renault Weyl groupoid construction with the Lawson-Lenz version of Exel's tight groupoid construction.  We do this by utilising only a weak algebraic fragment of the C*-algebra structure, namely its *-semigroup reduct.  Fundamental properties like local compactness are also shown to remain valid in general classes of *-rings.
\end{abstract}



\maketitle

\section{Introduction}

\subsection{Background}

Renault's groundbreaking thesis \cite{Renault1980} revealed the striking interplay between \'etale groupoids and the C*-algebras they give rise to.  Roughly speaking, the \'etale groupoid provides a more topological picture of the corresponding C*-algebra, with various key properties of the algebra, like nuclearity and simplicity (see \cite{AnantharamanDelaroucheRenault2000} and \cite{ClarkExelPardoSimsStarling2019}), being determined in a straightforward way from the underlying \'etale groupoid.  Naturally, this has led to the quest to find appropriate \'etale groupoid models for various C*-algebras.

Two general methods have emerged for finding such models, namely
\begin{enumerate}
\item Exel's tight groupoid construction from an inverse semigroup, and
\item Kumjian-Renault's Weyl groupoid construction from a Cartan C*-subalgebra
\end{enumerate}
(see \cite{Exel2008}, \cite{Kumjian1986} and \cite{Renault2008}).  However, both of these have their limitations.  For example, tight groupoids are always ample, which means the corresponding C*-algebras always have lots of projections.  On the other hand, the Weyl groupoid is always effective, which discounts many naturally arising groupoids.  Recently there has been a push to extend the Weyl groupoid construction in various directions, e.g. in \cite{CarlsenRuizSimsTomforde2017}, \cite{Resende2018}, \cite{ExelPitts2019} and \cite{KwasniewskiMeyer2019}, but even in these generalisations, some restrictions on the isotropy have remained.

Our goal is to unify these constructions in an elementary algebraic way which also eliminates these limitations.  The key is to utilise the *-semigroup structure of the *-normalisers of the Cartan C*-subalgebra just as one uses the inverse semigroup structure in the Lawson-Lenz approach to the tight groupoid construction (see \cite{LawsonLenz2013}).  In fact, one already sees *-semigroups playing a role in \cite{ExelPitts2019}.

More precisely, we show how to construct the Weyl groupoid via ultrafilters of *-normalisers with respect to appropriate transitive relations defined from the *-semigroup structure.  This idea of defining points from filters is a standard technique in point-free topology, going all the way back to the classic Stone and Wallman dualities (see \cite{Stone1936} and \cite{Wallman1938}).  In particular, this is the approach taken in \cite{Milgram1949} to recover a space $X$ from the semigroup structure of $C(X,\mathbb{R})$.  So from another perspective, what we are proposing is a non-commutative extension of \cite{Milgram1949}, which parallels our non-commutative extension of Stone duality in \cite{BiceStarling2018}.

\subsection{Motivation}

In defining the Weyl groupoid of a C*-algebra $A$, it is helpful to imagine that $A$ is already an algebra of continuous functions on some \'etale groupoid $G$ \textendash\, what we are trying to do is `reconstruct' or `recover' $G$ from the algebra $A$.  To motivate our method for doing this, it is instructive to first go back to the simplest case where the groupoid structure is trivial, i.e. where $G$ is just a space and $A$ is commutative.

\subsubsection{The Commutative Case}

Assume we have a locally compact Hausdorff space $X$ and consider the C*-algebra $C_0(X)$ of continuous $\mathbb{C}$-valued functions on $X$ vanishing at infinity, i.e.
\[C_0(X)=C_0(X,\mathbb{C})=\{f\in C(X):X\setminus f^{-1}[O]\textrm{ is compact, for every open }O\ni0\}.\]
The standard way to reconstruct $X$ from $C_0(X)$ is to identify points $x\in X$ with their evaluation functionals $\phi_x(f)=f(x)$ on $C_0(X)$.  One notes that every $\phi_x$ is a character and that every non-zero character is of the form $\phi_x$ for a unique $x\in X$.   Thus the character space $\Phi_{C_0(X)}$ with the usual weak* topology `recovers' $X$.

This also suggests that even an abstract commutative C*-algebra $A$ could be considered as an algebra of functions on its character space $\Phi_A$, which is exactly what the Gelfand representation theorem says.  Indeed, the standard Weyl groupoid associated to $A$ (with $A$ itself as the Cartan subalgebra) is exactly the character space of $A$.

Alternatively, we can identify points $x\in X$ with certain families of functions in $C_0(X)$.  One canonical choice would be the family of functions vanishing at $x$, i.e. the kernel of $\phi_x$.  Again one notes that these are maximal closed ideals and every maximal closed ideal is of this form, which allows us to recover $X$ as the space of maximal closed ideals in its Jacobson/hull-kernel topology.

Instead, we could identify each $x\in X$ with the family ${F}_x$ of functions with values in the unit interval $[0,1]$ taking the value $1$ on some neighbourhood of $x$, i.e.
\[{F}_x=\{f\in C_0(X,[0,1]):x\in\mathrm{int}(f^{-1}\{1\})\}.\]
Again one can immediately point out some basic properties of such families.  For example, if ${F}_x$ contains $f$ then ${F}_x$ also contains some $g$ dominated by $f$, by which we mean that $f$ is $1$ on the entire support of $g$,
\[\mathrm{supp}(g)=X\setminus g^{-1}\{0\}\subseteq f^{-1}\{1\}.\]
In algebraic terms, this can be expressed succinctly as $fg=g$.  Even for pairs $e,f\in{F}_x$, we can again find $g\in{F}_x$ dominated by both $e$ and $f$.  In the other direction, if $f\in{F}_x$ then every $g$ dominating $f$ again lies in ${F}_x$.  In other words, ${F}_x$ is a \emph{filter} with respect to the domination relation $\prec$ defined by
\[f\prec g\qquad\Leftrightarrow\qquad f=fg.\]
In fact each ${F}_x$ here is a maximal proper filter, i.e. an \emph{ultrafilter}, and again every ultrafilter has this form.  So the ultrafilters, again with a natural topology, provide an alternative reconstruction of the space $X$.

The important thing to note here is that, unlike maximal ideals or characters, these ultrafilters are determined from just the multiplicative semigroup structure.  This makes ultrafilters more suitable for reconstructing groupoids in a similar way.

\subsubsection{The Non-Commutative Case}\label{NonCommutativeCase}

Assume $G$ is a locally compact \'etale groupoid with Hausdorff unit space $G^0$ (see \autoref{*Semigroupoids} and \autoref{TopologicalGroupoids} below).  We call $O\subseteq G$ a \emph{bisection} if $OO^*\cup O^*O\subseteq G^0$ and let
\[S=\bigcup\{C_0(O,\mathbb{D}\setminus\{0\}):O\subseteq G\textrm{ is an open bisection}\}.\]
So each $a\in S$ is a function defined on an open bisection $\mathrm{dom}(a)\subseteq G$ with non-zero values in the complex disk $\mathbb{D}=\{\alpha\in\mathbb{C}:|\alpha|\leq1\}$ vanishing at infinity on $\mathrm{dom}(a)$.

First note that $S$ is a *-semigroup under the usual involution and product operations.  Specifically, for $g\in\mathrm{dom}(a)$ and $h\in\mathrm{dom}(b)$ with $g^*g=hh^*$,
\[a^*(g)=\overline{a(g^*)}\qquad\textrm{and}\qquad ab(gh)=a(g)b(h).\]
As $\mathrm{dom}(a)$ and $\mathrm{dom}(b)$ are bisections, $ab$ is a well defined function on another bisection $\mathrm{dom}(ab)=\mathrm{dom}(a)\mathrm{dom}(b)$, which is open because $G$ is \'etale.  As $a$ and $b$ are continuous, so is $ab$.  As the unit space is Hausdorff, products of compact subsets of $G$ are again compact so $ab$ vanishes at infinity on its domain, i.e. $ab\in S$.

The normal thing to do at this point would be to identify the partial functions in $S$ with the total functions taking $0$ values outside their original domain (even though the resulting total functions will not be continuous when $G$ itself is not Hausdorff).  One can then consider the *-algebra they generate under pointwise sums and the convolution product.  Completing this *-algebra in an appropriate norm, one then obtains a C*-algebra $A$.  Then $G$ can then be recovered from $A$ together with the Cartan subalgebra $C_0(G^0)$ via the Kumjian-Renault Weyl groupoid construction, as long as $G$ is Hausdorff, second countable and topologically principal.

However, we claim that to recover even non-Hausdorff $G$, regardless of its weight or isotropy, we just need the *-semigroup $S$ together with its `Cartan *-subsemigroup'
\[E=\{a\in S:\mathrm{dom}(a)\subseteq G^0\}.\]
To do this we proceed like in the commutative case.  Specifically, for any $a,b\in S$, we say that $b$ \emph{*-dominates} $a$, written $a\precsim b$, if $b$ takes values in the complex unit circle $\mathbb{T}=\{\alpha\in\mathbb{C}:|\alpha|=1\}$ on the entire domain of $a$, i.e.
\[a\precsim b\qquad\Leftrightarrow\qquad\mathrm{dom}(a)\subseteq b^{-1}[\mathbb{T}].\]
First we note $\precsim$ can be characterised by the *-semigroup structure of $S$ and $E$.

\begin{prp}\label{*DomChar}
For any $a,b\in S$,
\[a\precsim b\qquad\Leftrightarrow\qquad a=ab^*b\quad\textrm{and}\quad ab^*\in E.\]
\end{prp}

\begin{proof}
If $a\precsim b$ then, for any $g\in\mathrm{dom}(a)\subseteq b^{-1}[\mathbb{T}]$, $g\in\mathrm{dom}(b)$ and $b(g)\in\mathbb{T}$ so
\[a(g)=a(g)\overline{b(g)}b(g)=a(g)b^*(g^*)b(g)=ab^*b(gg^*g)=ab^*b(g),\]
i.e. $a=ab^*b$.  Moreover, $\mathrm{dom}(a)\subseteq\mathrm{dom}(b)$ implies
\[\mathrm{dom}(ab^*)=\mathrm{dom}(a)\mathrm{dom}(b)^*\subseteq\mathrm{dom}(b)\mathrm{dom}(b^*)\subseteq G^0,\]
i.e. $ab^*\in E$.

Conversely, assume $a=ab^*b$ and $ab^*\in E$ and take any $g\in\mathrm{dom}(a)$.  As $\mathrm{dom}(b^*b)=\mathrm{dom}(b)\mathrm{dom}(b)^*\subseteq G^0$, $a(g)=ab^*b(g)=a(h)b^*b(e)$, for some $e\in G^0$ with $g=he$, which can only mean $h=g$, $e=g^*g$ and $b^*b(e)=1$.  But $e=f^*f$, for some $f\in\mathrm{dom}(b)$, and then $gf^*\in\mathrm{dom}(a)\mathrm{dom}(b)^*=\mathrm{dom}(ab^*)\subseteq G^0$, as $ab^*\in E$.  Thus $g=f\in\mathrm{dom}(b)$ and $1=b^*b(e)=b^*(g^*)b(g)=\overline{b(g)}b(g)$, i.e. $b(g)\in\mathbb{T}$.  As $g$ was arbitrary, this shows that $a\precsim b$.
\end{proof}

Again consider \emph{ultrafilters}, i.e. maximal proper $T\subseteq S$ such that, for $a,b\in S$,
\[\tag{Filter}\label{FilterBeginning}a,b\in T\qquad\Leftrightarrow\qquad\exists c\in T\ (c\precsim a,b)\]
(note $\Rightarrow$ is saying $T$ is directed while $\Leftarrow$ is saying that $T$ is an up-set).  Let $\mathcal{U}(S)$ denote the space of ultrafilters with the topology generated by
\[\mathcal{U}_a=\{U\in\mathcal{U}(S):a\in U\},\]
for $a\in S$.  As each $U\in\mathcal{U}(S)$ is filter, $(\mathcal{U}_a)_{a\in S}$ is a basis for this topology.

We can now show that $\mathcal{U}(S)$ indeed recovers $G$ as claimed.

\begin{prp}\label{GroupoidRecovery}
We have a homeomorphism from $G$ onto $\mathcal{U}(S)$ given by
\[g\mapsto{S}_g=\{a\in S:g\in\mathrm{int}(a^{-1}[\mathbb{T}])\}.\]
This map is also a groupoid isomorphism, i.e. for all $g,h\in G$ with $g^*g=hh^*$,
\[S_{g^*}=(S_g)^*\qquad\textrm{and}\qquad{S}_{gh}=({S}_g{S}_h)^\precsim=\{a\succsim bc:b\in{S}_g\textrm{ and }c\in{S}_h\}.\]
\end{prp}

\begin{proof}
Take $g\in G$ and consider $S_g$.  If $a\succsim b\in S_g$ then $a\in S_g$, as
\[g\in\mathrm{int}(b^{-1}[\mathbb{T}])\subseteq\mathrm{dom}(b)\subseteq a^{-1}[\mathbb{T}].\]
If $a,b\in S_g$ then $g\in O=\mathrm{int}(a^{-1}[\mathbb{T}])\cap\mathrm{int}(b^{-1}[\mathbb{T}])$.  Urysohn's lemma yields $c\in S$ with $g\in\mathrm{int}(c^{-1}[1])$ and $\mathrm{dom}(c)\subseteq O$ which means $c\in S_g$ and $c\precsim a,b$.  Thus $S_g$ is a filter, which we claim is maximal.  If not, we would have some proper filter extension $F$, i.e. $S_g\subseteq F\neq S$.  For any $a\in F\setminus S_g$, $g\notin\mathrm{int}(a^{-1}[\mathbb{T}])$.  As $F$ is a filter, for any $b\in S_g(\neq\emptyset$, again by Urysohn's lemma$)$, we have $c\in F$ with $c\precsim a,b$ and hence $g\notin\mathrm{dom}(c)\supseteq c^{-1}[\mathbb{T}]$.  As $c^{-1}[\mathbb{T}]$ is compact subset of $\mathrm{dom}(b)$, Urysohn's lemma again yields $d\in S_g$ with $\mathrm{dom}(d)\cap c^{-1}[\mathbb{T}]=\emptyset$.  As $F$ is a filter, we have $e\in F$ with $e\precsim c,d$, but this can only mean $e$ is the empty function so $F=e^\precsim=S$, a contradiction.  Thus no such extension exists and hence $S_g$ is an ultrafilter.

On the other hand, assume $U$ is an ultrafilter and consider $C=\bigcap_{u\in U}\mathrm{dom}(u)$.  Fixing $a\in U$, we see that $C=\bigcap_{u\in U,u\precsim a}u^{-1}[\mathbb{T}]$ is a directed intersection of non-empty compact subsets of the Hausdorff subset $\mathrm{dom}(a)$ and must therefore contain some $g\in G$.  As $C=\bigcap_{u\in U}\mathrm{int}(u^{-1}[\mathbb{T}])$, it follows that $U\subseteq S_g$ and hence $U=S_g$, by maximality.  Thus $g\mapsto S_g$ is a bijection from $G$ onto $\mathcal{U}(S)$.  To see that it is a homeomorphism, note that Urysohn's lemma again shows that $(\mathrm{int}(a^{-1}[\mathbb{T}])_{a\in S}$ forms a basis for the topology of $G$ while, by definition, the sets
\[\mathcal{U}_a=\{U\in\mathcal{U}(S):a\in U\}=\{S_g:g\in\mathrm{int}(a^{-1}[\mathbb{T}])\}\]
form a basis for $\mathcal{U}(S)$.

For any $g\in G$, we immediately see that $S_{g^*}=(S_g)^*$ and $S_gS_h\subseteq S_{gh}$.  On the other hand, let $a\in S_{gh}$ and take $b\in S_{gh}$ with $b\precsim a$.  For any $c\in S_g$, $c^*b\in S_{g^*}S_{gh}\subseteq S_h$ so $cc^*b\in S_gS_h$ and
\[\mathrm{dom}(cc^*b)=\mathrm{dom}(c)\mathrm{dom}(c)^*\mathrm{dom}(b)\subseteq G^0\mathrm{dom}(b)=\mathrm{dom}(b)\subseteq a^{-1}[\mathbb{T}],\]
i.e. $cc^*b\precsim a$ and hence $a\in(S_gS_h)^\precsim$.  As $a$ was arbitrary, $S_{gh}=(S_gS_h)^\precsim$.
\end{proof}

\autoref{*DomChar} and \autoref{GroupoidRecovery} also suggest that, even starting from an abstract *-semigroup $S$ and *-subsemigroup $E$, we should be able to construct an \'etale groupoid $\mathcal{U}(S)$ of ultrafilters on $S$.  This is precisely what we will do.

\subsection{Outline}

First we set out some basic notation and terminology in \autoref{Preliminaries}.  In \autoref{Weyl*semigroup}, we introduce Weyl *-semigroups, the central objects under consideration throughout the paper.  We then give some motivating examples of Weyl *-semigroups in \autoref{Examples}.  In particular, \autoref{Closed*Subsemigroups} covers the case of commutative Cartan subalgebras in C*-algebras, while \autoref{ComplexAlgebras} even covers the non-commutative Cartan subalgebras from \cite{Exel2011} and \cite{KwasniewskiMeyer2019}.  We also note in \autoref{EquivalentNormalisers} that the *-normalisers of a commutative Cartan subalgebra $E$ coincide with the closure of the (non-*) normalisers of the self-adjoint part of $E$.

In \autoref{Relations} we introduce the relations that play a fundamental role in what follows and investigate their basic properties in relation to the *-semigroup structure.  We then examine cosets in \autoref{Cosets}, which provide the most natural general setting in which to define the \'etale groupoid structure.  In \autoref{Filters}, we restrict our attention to filters and eventually the ultrafilters which constitute our Weyl groupoid.  As a prelude to future work on C*-algebras and Fell bundles, in \autoref{Actions} we show how to construct the Weyl bundle in the general context of *-semigroups acting on uniform spaces, which then get represented as continuous sections of the bundle in \autoref{Representation}.

To say more about the Weyl groupoid, we need an ambient *-ring for our *-semigroup $S$.  In keeping with our algebraic approach, we first show in \autoref{Norms} how to define norm-like functions on any *-ring, which agree with the usual norm on C*-algebras.  This might seem like overkill, but we wish to illustrate that basic properties of the Weyl groupoid only depend on the *-ring structure.  Moreover, this should make it easy to apply our Weyl groupoid construction to algebraic cousins of C*-algebras like Leavitt path or Steinberg algebras, Baer *-rings, not to mention real C*-algebras.  Next, in \autoref{OrderStructure}, we examine elementary order properties of *-rings.  These allow us to show in \autoref{EtaleBasis} that the elements of $S$ naturally give rise to an inverse semigroup of open bisections in the Weyl groupoid.  Some additional lattice structure is introduced in \autoref{LatticeStructure} which we then utilise to show that the Weyl groupoid is locally compact in \autoref{WeylLocallyCompact}.

\subsection{Acknowledgements}

The author would like to thank Aidan Sims for some stimulating discussions on Weyl groupoids and also for making his \'etale groupoid notes \cite{Sims2017} available, as these were a big inspiration for the present paper.  Many thanks also go to Tatiana Shulman, Adam Skalski and the rest of the non-commutative geometry group at IMPAN for a very nice lecture series on \'etale groupoids in the first semester of 2019.  The author would also like to thank the anonymous referee who improved the paper with a number of helpful comments and suggestions.

\section{*-Semigroups}\label{StarSemigroups}

\subsection{Preliminaries}\label{Preliminaries}

First let us recall some basic definitions.

\begin{dfn}\label{*Semigroupoids}
A \emph{*-semigroupoid} is a set $S$ together with a total unary operation $a\mapsto a^*$ and a partial binary operation $(a,b)\mapsto ab$ such that, for all $a,b,c\in S$,
\begin{align}
\tag{Involution}a^{**}&=a.\\
\label{Associativity}\tag{Associativity}a(bc)&=(ab)c.\\
\label{Antihomomorphism}\tag{Antihomomorphism}(ab)^*&=b^*a^*.\\
\tag{Non-Trivial}aa^*&\textrm{ is defined}.
\end{align}

Call $u\in S$ a \emph{unit} if $ua=a$, whenever $ua$ is defined.  Units are denoted by
\[S^0=\{u\in S:u\textrm{ is a unit}\}.\]

Call $S$ a \emph{*-category} if every $a\in S$ has a range unit $\mathsf{r}(a)\in S$ such that
\[ab\textrm{ is defined}\qquad\Leftrightarrow\qquad a\mathsf{r}(b)\textrm{ is defined}.\]

Call $u\in S$ a \emph{unitary} if $uu^*$ and $u^*u$ are units.

A \emph{groupoid} is a *-category consisting entirely of unitaries.
\end{dfn}

Note here we interpret equations with partial operations like \eqref{Associativity} and \eqref{Antihomomorphism} as saying that `the left hand side is defined iff the right hand side is defined, in which case they are equal'.

\begin{rmk}\label{Semicategories}
Some authors use `semigroupoid' for what would more commonly be called a \emph{semicategory}, i.e. a subsemigroupoid of a category, where the composable pairs could be determined by an ambient collection of objects/units.  We will also be primarily concerned with *-semicategories, either groupoids or *-semigroups (which are *-semicategories with a single object).
\end{rmk}

First we note that units are not as one-sided as they may seem.

\begin{prp}
If $u$ is a unit then $u=u^*$ and $au=a$, for all $a\in S$.
\end{prp}

\begin{proof}
If $u$ is a unit then $ua=a$, for all $a\in S$.  In particular, $uu^*=u^*$ so
\[u=u^{**}=(uu^*)^*=u^{**}u^*=uu^*=u^*.\]
Thus if $au$ is defined then $(au)^*=u^*a^*=ua^*=a^*$ and hence $au=a$.
\end{proof}

From the range we also define source
\[\mathsf{s}(a)=\mathsf{r}(a^*).\]

\begin{prp}
If $C$ is a *-category then $a\mathsf{s}(a)$ is always defined and
\begin{align*}
ab\textrm{ is defined}\qquad&\Leftrightarrow\qquad\mathsf{s}(a)b\textrm{ is defined}.\\
abc\textrm{ is defined}\qquad&\Leftrightarrow\qquad \textrm{$ab$ and $bc$ are defined}.
\end{align*}
Moreover, $\mathsf{s}$ fixes units and $\mathsf{s}(a)$ is the unique unit such that $a\mathsf{s}(a)$ is defined.
\end{prp}

\begin{proof}
As $aa^*$ is always defined, so is $a\mathsf{r}(a^*)=a\mathsf{s}(a)$.  Next note $ab$ is defined iff $b^*a^*$ is defined iff $b^*\mathsf{r}(a^*)=b^*\mathsf{s}(a)$ is defined iff $\mathsf{s}(a)b=\mathsf{s}(a)^*b$ is defined, proving the first equivalence.

For the second, note that $abc$ being defined means that both $(ab)c$ and $a(bc)$ are defined so, in particular, both $ab$ and $bc$ are defined.  Conversely, assume $ab$ and $bc$ are defined.  Then $ab\mathsf{s}(ab)$ is defined and, in particular, $b\mathsf{s}(ab)$ is defined.  Thus $\mathsf{s}(b)\mathsf{s}(ab)$ is defined and, as sources are units, $\mathsf{s}(b)=\mathsf{s}(b)\mathsf{s}(ab)=\mathsf{s}(ab)$.  As $bc$ is defined, so is $\mathsf{s}(b)c=\mathsf{s}(ab)c$ and hence $abc$ is defined.

If $u$ is another unit such that $au$ is defined then $\mathsf{s}(a)u$ is defined and hence $\mathsf{s}(a)=\mathsf{s}(a)u=u$, proving uniqueness.  Also $u\mathsf{s}(u)$ is defined, for any unit $u$, and $u=u\mathsf{s}(u)=\mathsf{s}(u)$, so $\mathsf{s}$ fixes units.
\end{proof}

\begin{dfn}\label{TopologicalGroupoids}
A \emph{topological *-semigroupoid} is a *-semigroupoid $S$ under a topology making both the involution and product continuous (when defined).

We call $S$ \emph{open} if the product is also an open map.

We call $S$ \emph{\'etale} if the units $S^0$ are also open in $S$.
\end{dfn}

\'Etale groupoids are usually defined via local homeomorphisms, but this is just a special case of the above.

\begin{prp}\label{Resende}
A topological groupoid $G$ is \'etale iff the source $\mathsf{s}$ is an open map and hence a local homeomorphism(=open continuous locally injective).
\end{prp}

\begin{proof}
See \cite[Theorem 5.18]{Resende2007}.
\end{proof}

Let us call a subset $I$ of a *-category $C$ an \emph{ideal} if
\[\mathsf{r}(a)\in I\qquad\Leftrightarrow\qquad a\in I\qquad\Leftrightarrow\qquad\mathsf{s}(a)\in I.\]
Note any ideal is, in particular, a full *-subcategory.  Indeed, either of the equivalences above can be replaced by the condition that $I=I^*$.  Also, if $I$ is an ideal,
\[a\in I\quad\textrm{and}\quad ab\textrm{ is defined}\qquad\Rightarrow\qquad ab\in I,\]
as $\mathsf{r}(a)=\mathsf{r}(ab)$.  If $C$ is a groupoid and $I=I^*$, this even characterises the ideals, as then $a\in I$ implies $\mathsf{r}(a)=aa^*\in I$, which in turn implies $a=\mathsf{r}(a)a\in I$.

\begin{prp}\label{OpenIdeals}
If $C$ is an \'etale or open *-category then so are its ideals $I$.
\end{prp}

\begin{proof}
If the product and involution are continuous on $C$ then so are their restrictions to $I$, and if $C^0$ is open in $C$ then $I^0=I\cap C^0$ is open in $I$.  Now assume the product is an open map, i.e. $ON$ is open, for any open $O,N\subseteq C$.  If $I$ is a *-subcategory then $(O\cap I)(N\cap I)\subseteq ON\cap I$.  If $I$ is an ideal then the reverse inclusion also holds, as $ab\in I$ implies $r(a)=r(ab)\in I$ and $s(b)=s(ab)\in I$ and hence $a,b\in I$.  Thus $(O\cap I)(N\cap I)=ON\cap I$ is open in $I$, so the product on $I$ is again an open map.
\end{proof}

Given a subset $T$ of a *-semigroupoid $S$, we define
\begin{align}
\tag{*-Squares}|T|^2&=\{a^*a:a\in T\}.\\
\tag{Self-Adjoints}T_\mathrm{sa}&=\{a\in T:a=a^*\}.\\
\tag{Partial Isometries}T_\mathrm{pi}&=\{a\in T:a=aa^*a\}.\\
\tag{*-Normalisers}T_*&=\{a\in S:a^*Ta\cup aTa^*\subseteq T\}.\\
\tag{Normalisers}T_\bullet&=\{a\in S:aT=Ta\}.\\
\tag{Centralisers}T'&=\{a\in S:\forall t\in T\ (at=ta)\}=\bigcap_{t\in T}t_\bullet.
\end{align}

A \emph{*-semigroup} is a *-semigroupoid where the product is defined on all pairs.

A subset $T$ of a *-semigroup $S$ is a \emph{*-subsemigroup} if $TT\subseteq T=T^*$.

A subset $T$ of a *-semigroup $S$ is \emph{*-normal} if $T_*=S$, i.e. $a^*Ta\subseteq T$, for $a\in S$.

\

We now define the *-semigroups that are central to the present paper.

\begin{dfn}\label{Weyl*semigroup}
A \emph{Weyl *-semigroup} $(S,E)$ is a *-semigroup $S$ with a distinguished *-normal *-subsemigroup $E$ whose centre contains the *-squares, i.e.
\[\tag{Central *-Squares}|S|^2\subseteq E\cap E'.\]
\end{dfn}

\subsection{Examples}\label{Examples}

First we give some motivating examples of Weyl *-semigroups.

\begin{xpl}\label{Groups}
If $S$ is a group (with $s^*=s^{-1}$) and $E$ is any normal subgroup then $(S,E)$ is a Weyl *-semigroup, as $|S|^2=\{e\}$ is the unit of $S$, which is certainly central in $S$ and hence in $E$.
\end{xpl}

\begin{xpl}\label{InverseSemigroups}
If $S$ is an inverse semigroup then the idempotents $|S|^2$ form a commutative *-normal *-subsemigroup so $(S,|S|^2)$ is a Weyl *-semigroup.
\end{xpl}

In fact, every Weyl *-semigroup contains a canonical inverse semigroup, namely its *-subsemigroup of partial isometries $S_\mathrm{pi}$.  To see this, note that if $a\in S_\mathrm{pi}$ then $a=aa^*a$ and hence $a^*=(aa^*a)^*=a^*aa^*$ so $a^*\in S_\mathrm{pi}$.  Similarly, if $a,b\in S_\mathrm{pi}$ then, as $|S|^2$ is commutative, $ab=aa^*abb^*b=abb^*a^*ab=ab(ab)^*ab$, so $ab\in S_\mathrm{pi}$.  Moreover, if $a\in S_\mathrm{pi}$ then $a^*a=a^*aa^*a$, i.e. every element of $|S_\mathrm{pi}|^2$ is idempotent.  Conversely, if $a\in S$ is idempotent then $a=aa$ immediately yields $a^*=(aa)^*=a^*a^*$ and hence $a=aa^*a=aa^*a^*a=a^*aaa^*=a^*aa^*=a^*$, as $|S|^2$ is commutative, so $a=aa=a^*a\in|S_\mathrm{pi}|^2$.  As the idempotents $|S_\mathrm{pi}|^2$ commute, $S_\mathrm{pi}$ is an inverse semigroup.  In particular, inverse semigroups are precisely the Weyl *-semigroups consisting entirely of partial isometries.

\begin{xpl}\label{CommutativeSemigroups}
If $S$ is a commutative semigroup then we can define $s^*=s$, for all $s\in S$.  For any subsemigroup $E$ containing all the (*-)squares $S^2=|S|^2$, we see that $(S,E)$ is then a Weyl *-semigroup. 
\end{xpl}

\begin{xpl}\label{Commutative*Subsemigroups}
Assume $A$ is a *-semigroup with a commutative *-subsemigroup $E$.  Taking $E'$ and $E_*$ within $A$, we immediately see that $E\subseteq E'$ and $E\subseteq E_*$.  Also, $E_*$ is a *-subsemigroup of $A$ and $E$ is certainly *-normal in $E_*$.  If $E$ contains a unit for $A$ then $|E_*|^2\subseteq E$ so $(E_*,E)$ is a Weyl *-semigroup.

When $A$ is also commutative, $E_*=\{a\in A:a^*a,aa^*\in E\}$.  For example, we could take $A=\mathbb{C}^2$ (with the pointwise product) and
\[E=\{(\alpha,\alpha):\alpha\in\mathbb{C}\}.\]
If we also take $(\alpha,\beta)^*=(\alpha,\beta)$, as suggested in \autoref{CommutativeSemigroups}, then
\[E_*=E\cup\{(\alpha,-\alpha):\alpha\in\mathbb{C}\}.\]
On the other hand, if we consider $A$ as a C*-algebra, i.e. $(\alpha,\beta)^*=(\overline\alpha,\overline\beta)$, then
\[E_*=\{(\alpha,t\alpha):\alpha\in\mathbb{C}\textrm{ and }t\in\mathbb{T}\},\]
Similar remarks apply when we consider the C*-algebra $A=C([0,1])$ and C*-subalgebra $E=\{f\in A:f(0)=f(1)\}$, as in \cite[Example 22]{ExelPitts2019}.
\end{xpl}

\begin{xpl}\label{Closed*Subsemigroups}
Assume $A$ is a topological *-semigroup (so the product and involution are continuous) and $E$ is a closed commutative *-subsemigroup.  If $E$ contains a left or right approximate unit for $A$, i.e. a net $(e_\lambda)\subseteq E$ with $ae_\lambda\rightarrow a$ or $e_\lambda a\rightarrow a$, for all $a\in A$, again $(E_*,E)$ is a Weyl *-semigroup.
\end{xpl}

\begin{xpl}\label{ConditionalExpectations}
Assume $\Phi$ is a *-preserving $E$-equivariant map on a *-semigroup $A$, where $E=\Phi[A]$ is commutative.  In other words, for all $a\in A$ and $e\in E$,
\[\Phi(a^*)=\Phi(a)^*\qquad\textrm{and}\qquad\Phi(ea)=e\Phi(a),\]
so $E=\Phi[A]$ is a *-subsemigroup of $A$ and, for all $a\in A$ and $e\in E$,
\[\Phi(ae)=\Phi(e^*a^*)^*=(e^*\Phi(a^*))^*=\Phi(a)e.\]
Consider the *-normalisers of $\Phi$,
\[\Phi_*=\{s\in A:\forall a\in A\ (\Phi(s^*as)=s^*\Phi(a)s\textrm{ and }\Phi(sas^*)=s\Phi(a)s^*\}.\]
Note $\Phi_*$ is *-subsemigroup of $A$ containing $E$, as $\Phi(eae^*)=e\Phi(ae^*)=e\Phi(a)e^*$, for any $e\in E$.   Also, for any $e\in E$, we have $a\in A$ with $e=\Phi(a)$ and hence 
$s^*es=s^*\Phi(a)s=\Phi(s^*as)\in E$, for any $s\in\Phi_*$, showing that $E$ is *-normal in $\Phi_*$.  If $A$ is unital and $\Phi(1)=1$ then, for any $s\in\Phi_*$, we see that $s^*s=s^*1s=s^*\Phi(1)s=\Phi(s^*1s)=\Phi(s^*s)\in E$ so $|\Phi_*|^2\subseteq E$.  Alternatively, if $A$ is a Hausdorff topological *-semigroup, $\Phi$ is continuous and idempotent and $E$ contains a left or right approximate unit $(e_\lambda)$ for $A$ then
\[s^*s\leftarrow s^*e_\lambda s=s^*\Phi(e_\lambda)s=\Phi(s^*e_\lambda s)\rightarrow\Phi(s^*s)\in E,\]
so again $|\Phi_*|^2\subseteq E$, i.e. in either case $(\Phi_*,E)$ is a Weyl *-semigroup.
\end{xpl}

In particular, \autoref{Closed*Subsemigroups} and \autoref{ConditionalExpectations} apply when $A$ is a C*-algebra with commutative Cartan subalgebra $E=\Phi[A]$ (see \cite[Definition 5.1]{Renault2008}).  In this case, we know that $\Phi_*=E_*$ by going through Renault's Weyl groupoid construction, although it would be nice to have a more direct proof of this.  For more general conditional expectations onto non-maximal commutative C*-subalgebras, we only have $\Phi_*\subseteq E_*$.  For example, if $\Phi(a)=\phi(a)1$, for some pure state $\phi(a)=\langle av,v\rangle$ on $M_2$ then $E=\mathbb{C}1$ and $E_*$ consists of all scalar multiples of all unitaries $u\in M_2$, while to be in $\Phi_*$, $u$ must also fix the $1$-dimensional subspace spanned by $v$.

\begin{xpl}\label{ComplexAlgebras}
Assume $A$ is a complex *-algebra with a *-subalgebra $E$ and
\[E'\subseteq E.\]
Then $(E_{\mathrm{sa}\bullet},E)$ is a Weyl *-semigroup (recall $E_{\mathrm{sa}\bullet}$ consists normalisers of the self-adjoint elements of $E)$.  Indeed, for any $a\in E_{\mathrm{sa}\bullet}$ and $e\in E_\mathrm{sa}$, we have $f\in E_\mathrm{sa}$ with $ae=fa$ and hence $ea^*=a^*f$ so $a^*ae=a^*fa=ea^*a$, i.e.
\[|E_{\mathrm{sa}\bullet}|^2\subseteq(E_\mathrm{sa})'=(E_\mathrm{sa}+iE_\mathrm{sa})'=E'\subseteq E.\]
For any $a\in E_{\mathrm{sa}\bullet}$, it follows that $a^*E_\mathrm{sa}a=a^*aE_\mathrm{sa}\subseteq|E_{\mathrm{sa}\bullet}|^2E_\mathrm{sa}\subseteq EE\subseteq E$
so $a^*Ea=a^*(E_\mathrm{sa}+iE_\mathrm{sa})a\subseteq E$.  Thus $E_{\mathrm{sa}\bullet}\subseteq E_*$, i.e. $E$ is *-normal in $E_{\mathrm{sa}\bullet}$.

If $A$ is also topological then $E$ and $E_*$ are automatically closed.  The same can not be said for $E_{\mathrm{sa}\bullet}$, although we can always take $S=\overline{E_{\mathrm{sa}\bullet}}$, as $\overline{E_{\mathrm{sa}\bullet}}\subseteq\overline{E_*}=E_*$ and $|\overline{E_{\mathrm{sa}\bullet}}|^2\subseteq\overline{|E_{\mathrm{sa}\bullet}|^2}\subseteq\overline{E}=E$.
\end{xpl}

Again, \autoref{ComplexAlgebras} applies when $A$ is a C*-algebra with Cartan subalgebra $E$, even when $E$ is a non-commutative Cartan subalgebra in the sense of \cite{Exel2011}.  In C*-algebras, the focus is mainly on *-normalisers rather than normalisers, but these turn out to be the same thing for commutative Cartan subalgebras.

\begin{thm}\label{EquivalentNormalisers}
If $A$ is a C*-algebra with a commutative C*-subalgebra $E$ containing an approximate unit for $A$ then $E_*\subseteq\overline{E_{\mathrm{sa}\bullet}}$.\\
If $E$ is also maximal commutative then $E_*=\overline{E_{\mathrm{sa}\bullet}}$.
\end{thm}

\begin{proof}
Take $a\in E_*$.  As $E$ contains an approximate unit $(a_\lambda)$ for $A$,
\[a^*a=\lim_\lambda a^*a_\lambda a\in\overline{E}=E.\]
Representing $A$ on a Hilbert space $H$, we obtain a polar decomposition
\[a=u\sqrt{a^*a}=\sqrt{aa^*}u,\]
for some $u\in\mathcal{B}(H)$ with $uu^*a=a=au^*u$.  For any $\epsilon>0$, we have continuous $f$ and $g$ on $\mathbb{R}_+$ fixing $0$ such that $|f(r)-\sqrt{r}|<\epsilon$, for all $r\in\mathbb{R}_+$, and $fg=f$.  As $g(a^*a)=\lim_na^*ap_n(a^*a)$, for some sequence of polynomials $p_n$, we can set
\[c=g(aa^*)u=ug(a^*a)=\lim_nua^*ap_n(a^*a)=\lim_na\sqrt{a^*a}p_n(a^*a)\in aE\in E_*E\subseteq E_*.\]
Also set $b=f(aa^*)u=uf(a^*a)\in A$ so $cc^*b=g^2(aa^*)f(aa^*)u=f(aa^*)u=b$ and
\[c^*b=u^*g(aa^*)f(aa^*)u=u^*f(aa^*)u=\sqrt{u^*f^2(aa^*)u}=\sqrt{b^*b}=f(a^*a)\in E.\]
As $E$ is commutative, for any $e\in E_\mathrm{sa}$, $c^*be=ec^*b$ and hence, as $c\in E_*$,
\[be=cc^*be=cec^*b\in E_\mathrm{sa}b.\]
Thus $bE_\mathrm{sa}\subseteq E_\mathrm{sa}b$ and, likewise, $E_\mathrm{sa}b\subseteq bE_\mathrm{sa}$, i.e. $b\in E_{\mathrm{sa}\bullet}$.  As $\epsilon>0$ was arbitrary, we can take a sequence $(b_n)\subseteq E_{\mathrm{sa}\bullet}$ with $b_n\rightarrow a$ so $E_*\subseteq\overline{E_{\mathrm{sa}\bullet}}$.

If $E$ is maximal commutative then $E'\subseteq E$ so $\overline{E_{\mathrm{sa}\bullet}}\subseteq E_*$, as in \autoref{ComplexAlgebras}.
\end{proof}

We note that observations similar to those in \autoref{ComplexAlgebras} and \autoref{EquivalentNormalisers} can be found in \cite[Proposition 3.3 and 3.4]{DonsigPitts2008}.

The primary motivating examples for the Weyl groupoid construction come from continuous functions on an \'etale groupoid $G$.  In \autoref{NonCommutativeCase}, we considered functions taking values in the complex unit disk, but we could also consider values in a more general *-semigroup, or even sections of a *-category bundle over $G$.

\begin{xpl}\label{GroupoidBundles}
Again take a locally compact \'etale groupoid $G$ with Hausdorff unit space $G^0$.  Now assume we have an additional topological *-category $F$ together with a continuous \emph{*-isocofibration} $\pi:F\rightarrow G$, i.e. a functor preserving the involution which is also injective on the units $F^0$ so, for all $e,f\in F$,
\[\mathsf{s}(e)=\mathsf{r}(f)\qquad\Leftrightarrow\qquad\mathsf{s}(\pi(e))=\mathsf{r}(\pi(f)),\]
 in which case $\pi(ef)=\pi(e)\pi(f)$.  Further assume that (like with groupoids)
 \begin{equation}\label{eefAss}
 \mathsf{s}(e)=\mathsf{s}(f)=\mathsf{r}(f)\qquad\Rightarrow\qquad e^*ef=fe^*e.
 \end{equation}
 Assume we also have a seminorm $\|\cdot\|:F\rightarrow\mathbb{R}_+$, i.e. for all $e,f\in F$ with $\mathsf{s}(e)=\mathsf{r}(f)$,
\[\|f\|=\|f^*\|\qquad\textrm{and}\qquad\|ef\|\leq\|e\|\|f\|.\]
Let $S$ consist of (partial) continuous sections on the open bisections $\mathcal{B}(G)$ of $G$,
\[S=\{a\in C(O,F):O\in\mathcal{B}(G),\|a\|\in C_0(O)\textrm{ and }\pi(a(g))=g,\textrm{ for all }g\in O\}.\]
As in \autoref{NonCommutativeCase}, one immediately verifies that $S$ is a *-semigroup with *-subsemigroup
\[E=\{a\in S:\mathrm{dom}(a)\subseteq G^0\}.\]
For any $a\in E$ and $b\in S$, we see that
\[\mathrm{dom}(b^*ab)=\mathrm{dom}(b)^*\mathrm{dom}(a)\mathrm{dom}(b)\subseteq\mathrm{dom}(b)^*G^0\mathrm{dom}(b)\subseteq\mathrm{dom}(b)^*\mathrm{dom}(b)\subseteq G^0,\]
so $b^*ab\in E$ and hence $E$ is *-normal.  Also, thanks to \eqref{eefAss},
\[|S|^2\subseteq\{a\in E:\mathrm{ran}(a)\subseteq|F|^2\}\subseteq E',\]
i.e. elements of $|S|^2$ commute with those of $E$, so $(S,E)$ is a Weyl *-semigroup.
\end{xpl}

To see that \autoref{GroupoidBundles} encompasses \autoref{NonCommutativeCase}, just take
\[F=G\times(\mathbb{D}\setminus\{0\})\]
with its canonical groupoid structure $(g,\alpha)(h,\beta)=(gh,\alpha\beta)$ (when $\mathsf{s}(g)=\mathsf{r}(h)$), norm $\|(g,\alpha)\|=|\alpha|$ and $\pi(g,\alpha)=g$.  \autoref{GroupoidBundles} also subsumes \autoref{Groups} \textendash\, given a discrete group $F$ with normal subgroup $N$, we can take $G=F/N$, let $\pi:F\rightarrow G$ be the canonical quotient map and define $\|f\|=1$, for all $f\in F$.  Then each partial section $a\in S$ just selects a single element from a single coset, the elements of $E$ being those which select an element of $N$, i.e. we get back $F$ and $N$ when we identify partial sections with their images.

\subsection{Relations}\label{Relations}

To avoid unnecessary repetition, let us assume throughout the rest of the paper that $(S,E)$ is a given Weyl *-semigroup as per \autoref{Weyl*semigroup}.

\begin{ass}\label{SemigroupAssumption}
\textbf{$(S,E)$ is a Weyl *-semigroup.}
\end{ass}

\begin{dfn}
We consider the following relations on $S$.
\begin{align}
\label{Domination}\tag{Domination}a\prec b\qquad&\Leftrightarrow\qquad a=ab.\\
\label{Compatibility}\tag{Compatibility}a\sim b\qquad&\Leftrightarrow\qquad ab^*\in E.\\
\label{NaturalOrder}\tag{*-Domination}a\precsim b\qquad&\Leftrightarrow\qquad b\sim a\prec b^*b.
\end{align}
\end{dfn}

\begin{rmk}
The domination relation was considered on commutative semigroups in \cite{Milgram1949} and on the positive cone of C*-algebras in \cite[Definition II.3.4.3]{Blackadar2017}.  Our $\sim$ agrees with left compatibility, as defined in \cite[\S1.4]{Lawson1998} (right compatibility is $a^*\sim b^*$), while *-domination also agrees with the natural order defined in \cite[\S1.4]{Lawson1998} for inverse semigroups.  But unlike in inverse semigroups, *-domination is not always a partial order, it is only a transitive relation.
\end{rmk}

\begin{xpl}
If we consider the C*-algebra $S=E=C([0,1])$ then $\sim$ is trivial, i.e. $a\sim b$, for all $a,b\in S=E$, while $\prec$ and $\precsim$ can be characterised as
\begin{align*}
a\prec b\qquad&\Leftrightarrow\qquad\mathrm{supp}(a)\subseteq b^{-1}\{1\}.\\
a\precsim b\qquad&\Leftrightarrow\qquad\mathrm{supp}(a)\subseteq b^{-1}[\mathbb{T}].
\end{align*}
In particular, $\precsim$ is not reflexive, e.g. $\mathrm{id}\not\precsim\mathrm{id}$ where $\mathrm{id}(x)=x$, for all $x\in[0,1]$.
\end{xpl}

\begin{xpl}
When $S$ is a group with identity $e$ and normal subgroup $E$,
\begin{align*}
a\prec b\qquad&\Leftrightarrow\qquad e=b.\\
a\sim b\qquad&\Leftrightarrow\qquad a\precsim b\\
&\Leftrightarrow\qquad a\textrm{ and $b$ are in the same $E$-coset}.
\end{align*}
So $\precsim$ here is not antisymmetric unless $E$ is the trivial one-element subgroup.
\end{xpl}

\begin{prp}
For all $a,b,c\in S$,
\begin{align}
\label{absimc}\tag{$\sim$-Auxiliary}a\precsim b\sim c\qquad&\Rightarrow\qquad a\sim c.\\
\label{abprecc}\tag{$\prec$-Auxiliary}a\precsim b\prec c\qquad&\Rightarrow\qquad a\prec c.\\
\label{Transitivity}\tag{Transitivity}a\precsim b\precsim c\qquad&\Rightarrow\qquad a\precsim c.
\end{align}
\end{prp}

\begin{proof}\
\begin{itemize}
\item[\eqref{absimc}] If $a\precsim b\sim c$ then $ac^*=ab^*bc^*\in EE\subseteq E$, i.e. $a\sim c$.
\item[\eqref{abprecc}] If $a\precsim b\prec c$ then $ac=ab^*bc=ab^*b=a$, i.e. $a\prec c$.
\item[\eqref{Transitivity}] If $a\precsim b\precsim c$ then $a\precsim b\sim c$ and $a\precsim b\prec c^*c$ so $c\sim a\prec c^*c$, i.e. $a\precsim c$.
\end{itemize}
\end{proof}

\begin{prp}
For all $a,b\in S$,
\[\label{*Invariance}\tag{*-Invariance}a\precsim b\qquad\Leftrightarrow\qquad a^*\precsim b^*.\]
\end{prp}

\begin{proof}
Assume $a\precsim b$ so $ab^*\in E$ and $a=ab^*b$ or, equivalently, $ba^*\in E$ and $a^*=b^*ba^*$.  Thus $a^*b=b^*ba^*b\in b^*Eb\subseteq E$, i.e. $a^*\sim b^*$, and $a^*bb^*=b^*ba^*bb^*=b^*bb^*ba^*=a^*$, as $ba^*\in E$ and $bb^*\in |S|^2$ commute, i.e. $a^*\prec bb^*$ so $a^*\precsim b^*$.
\end{proof}

\begin{prp}
*-Domination preserves the product, i.e. for $a,b,c,d\in S$,
\[\tag{Product Preserving}\label{ProductPreserving}a\precsim b\quad\textrm{and}\quad c\precsim d\qquad\Rightarrow\qquad ac\precsim bd.\]
\end{prp}

\begin{proof}
If $a\precsim b$ and $c\precsim d$ then $cd^*\in E$ commutes with $b^*b\in|S|^2$ and hence $ac(bd)^*bd=acd^*b^*bd=ab^*bcd^*d=ac$, i.e. $ac\prec(bd)(bd)^*$.  As $E$ is also *-normal,
\[ac(bd)^*=acd^*b^*= ab^*bcd^*b^*\in EbEb^*\subseteq EE\subseteq E.\]
\end{proof}

\begin{prp}
For all $a,b,c\in S$,
\begin{equation}\label{TripleProduct}
a\precsim b\quad\textrm{and}\quad a\prec cc^*\qquad\Rightarrow\qquad a\precsim bcc^*\quad\textrm{and}\quad ac\precsim bc.
\end{equation}
\end{prp}

\begin{proof}
If $cc^*\succ a\precsim b$ then $a(bcc^*)^*=acc^*b^*=ab^*\in E$, i.e. $a\sim bcc^*$, and $a(bcc^*)^*bcc^*=acc^*b^*bcc^*=a$, so $a\precsim bcc^*$.  Likewise, $ac(bc)^*=acc^*b^*\in E$, i.e. $ac\sim bc$ and $ac(bc)^*bc=acc^*b^*bc=ac$, i.e. $ac\prec(bc)^*bc$ and hence $ac\precsim bc$.
\end{proof}

\begin{prp}
For any $a\in S$,
\begin{align}
\label{EDownClosed}a\precsim e\in E\qquad&\Rightarrow\qquad a\in E.\\
\label{Esucc}a^*a\succ e\in E\qquad&\Rightarrow\qquad e\precsim a^*a.\\
\label{Esuccsim}a\succsim e\in E\qquad&\Rightarrow\qquad e\prec aa^*.
\end{align}
\end{prp}

\begin{proof}\
\begin{itemize}
\item[\eqref{EDownClosed}] If $a\precsim e\in E$ then $ae^*\in E$ and $a=ae^*e\in EE\subseteq E$.
\item[\eqref{Esucc}] If $a^*a\succ e\in E$ then $e(a^*a)^*a^*a=ea^*aa^*a=e$, i.e. $e\prec(a^*a)^*a^*a$, and $e(a^*a)^*=ea^*a=e\in E$, i.e. $e\sim a^*a$, so $e\precsim a^*a$.
\item[\eqref{Esuccsim}] If $a\succsim e\in E$ then $e^*\precsim a^*$ so, as $|S|^2\subseteq E'$, $eaa^*=aa^*e=e$.
\end{itemize}
\end{proof}

\begin{prp}
For all $a,b,c\in S$ and $e\in E$,
\begin{align}
\label{EBelow}a\precsim b\qquad&\Rightarrow\qquad ea\precsim b\quad\textrm{and}\quad ae\precsim b.\\
\label{SBelow}a\prec b\qquad&\Rightarrow\qquad ca\prec b.
\end{align}
\end{prp}

\begin{proof}\
\begin{itemize}
\item[\eqref{EBelow}] Note $eab^*\in EE\subseteq E$, i.e. $ea\sim b$, and $eab^*b=ea$, i.e. $ea\prec b^*b$, so $ea\precsim b$.  Then $a\precsim b$ implies $a^*\precsim b^*$ and hence $e^*a^*\precsim b^*$, i.e. $ae\precsim b$.

\item[\eqref{SBelow}] Just note $a=ab$ implies $ca=cab$.
\end{itemize}
\end{proof}

\subsection{Cosets}\label{Cosets}

For any binary relation $\sqsubset$ on $S$ and $T\subseteq S$, let
\[T^\sqsubset=\{a\in S:\exists t\in T\ (t\sqsubset a)\}.\]
For example, \eqref{EDownClosed} could be rewritten as $E^\succsim\subseteq E$, i.e. $E$ is a down-set w.r.t. $\precsim$.

\begin{dfn}[{\cite[\S1.4]{Lawson1998}}]
We call $C\subseteq S$ an \emph{atlas} or \emph{coset} if
\begin{align*}
\label{Atlas}\tag{Atlas}CC^*C\cup C&\subseteq C^\precsim.\\
\label{Coset}\tag{Coset}CC^*C\subseteq C&=C^\precsim.
\end{align*}
\end{dfn}

Note $C$ is an atlas iff $C^\precsim$ is a coset containing $C$, i.e. atlases are coinitial subsets of cosets.  Sometimes even weak versions of \eqref{Coset} suffice, as in the following.

\begin{prp}\label{*Directed}
If $C|C|^2\subseteq C\subseteq C^\precsim$ and $a,b\in C$ then we have $c,d\in C$ with
\[c\precsim a,\quad d\precsim b,\quad c\prec b^*b,\quad d\prec a^*a\quad\textrm{and}\quad c^*c=d^*d.\]
\end{prp}

\begin{proof}
As $C\subseteq C^\precsim$, we have $e,f,g,h\in C$ with $e\precsim f\precsim a$ and $g\precsim h\precsim b$.  Let $c=fe^*eg^*g\in C|C|^2|C|^2\subseteq C$ and $d=hg^*ge^*e\in C|C|^2|C|^2\subseteq C$.  Note
\[c^*c=g^*ge^*ef^*fe^*eg^*g=g^*ge^*ee^*eg^*g=e^*eg^*gg^*ge^*e=e^*eg^*gh^*hg^*ge^*e=d^*d.\]
By \eqref{EBelow}, $f\precsim a$ and $c\in fE$ yields $c\precsim a$, while $h\precsim b$ and $d\in hE$ yields $d\precsim b$.  Also $c\prec b^*b$, as $g\prec b^*b$ and $c\in Sg$, and $d\prec a^*a$, as $e\prec a^*a$ and $d\in Se$.
\end{proof}

\begin{prp}\label{CaAtlas}
If $\{aa^*\}^\succ\cap(C^*C)^\precsim\neq\emptyset$ and $C$ is an atlas then so are $Ca$, $Caa^*$\nolinebreak and
\begin{equation}\label{Caa*}
C^\precsim=(Caa^*)^\precsim.
\end{equation}
\end{prp}

\begin{proof}
Take $b\in\{aa^*\}^\succ\cap(C^*C)^\precsim$.  For any element $c$ of the coset $C^\precsim$,
\[cb\in C^\precsim(C^*C)^\precsim\subseteq(CC^*C)^\precsim\subseteq C^{\precsim\precsim}\subseteq C^\precsim.\]
Then \autoref{*Directed} yields $f\in C^\precsim$ with $f\precsim c$ and $f\prec b^*c^*cb$.  Thus we have $g\in C$ with $g\precsim f$ and hence $g=gf^*f=gf^*fb^*c^*cb=gf^*fb^*c^*cbaa^*=gaa^*$,
i.e. $g\prec aa^*$.  By \eqref{TripleProduct}, $g\precsim caa^*$ and $ga\precsim ca$.  As $c$ was arbitrary, $C^\precsim aa^*\subseteq C^\precsim$ so
\[Caa^*\subseteq C^\precsim aa^*\subseteq C^\precsim\subseteq(Caa^*)^\precsim\subseteq C^{\precsim\precsim}\subseteq C^\precsim,\]
 i.e. $Caa^*$ is an atlas with $(Caa^*)^\precsim=C^\precsim$.  Again as $c$ was arbitrary, $C^\precsim a\subseteq(Ca)^\precsim$ so $Ca\subseteq C^\precsim a\subseteq(Ca)^\precsim$ and $Caa^*C^*Ca\subseteq C^\precsim C^{\precsim*}C^\precsim a\subseteq C^\precsim a\subseteq(Ca)^\precsim$, showing that $Ca$ is also an atlas.
\end{proof}

\begin{prp}
Take any $T,U,C\subseteq S$.
\begin{align}
\label{TUC<=TU}&\textrm{If $C$ is an atlas,}&TC&\subseteq C^\precsim&&\Leftarrow\qquad T\subseteq(CC^*)^\precsim.\\
\label{TUC=>TU}&\textrm{If $C\neq\emptyset$,}\quad&TC&\subseteq C&&\Rightarrow\qquad T^\precsim\subseteq(CC^*)^\precsim.\\
\label{TUC=>TC}&\textrm{If }U\neq\emptyset,&U^*TC&\subseteq C&&\Rightarrow\qquad(TC)^\precsim\subseteq(UC)^\precsim.\\
\label{TUC=>Proper}&\textrm{If }U^\precsim\neq\emptyset\textrm{ and }C^\precsim\neq S,&UTC&\subseteq C&&\Rightarrow\qquad(TC)^\precsim\neq S.
\end{align}
\end{prp}

\begin{proof}\
\begin{itemize}
\item[\eqref{TUC<=TU}] Note $TC\subseteq(CC^*)^\precsim C^\precsim\subseteq(CC^*C)^\precsim\subseteq C^{\precsim\precsim}\subseteq C^\precsim$, as $C$ is an atlas.

\item[\eqref{TUC=>TU}] By \eqref{EBelow} and $|S|^2\subseteq E$, $TC\subseteq C\neq\emptyset$ implies $T^\precsim\subseteq(TCC^*)^\precsim\subseteq(CC^*)^\precsim$.

\item[\eqref{TUC=>TC}] Again by \eqref{EBelow} and $|S|^2\subseteq E$, $(TC)^\precsim\subseteq(UU^*TC)^\precsim\subseteq(UC)^\precsim$.

\item[\eqref{TUC=>Proper}] Take $a\in S\setminus C^\precsim$ and $c\in U^\precsim$.  If $(TC)^\precsim=S$ then we would have $b\precsim a$ so
\[a\succsim cc^*b\in U^\precsim S=U^\precsim(TC)^\precsim\subseteq(UTC)^\precsim\subseteq C^\precsim,\]
 by \eqref{EBelow}, so $a\in C^{\precsim\precsim}\subseteq C^\precsim$, a contradiction.
\end{itemize}
\end{proof}

Let us denote the non-empty cosets in $S$ by
\[\mathcal{C}(S)=\{C\subseteq S:CC^*C\subseteq C\subseteq C^\precsim\neq\emptyset\}.\]

\begin{thm}
The cosets $\mathcal{C}(S)$ form a groupoid under the inverse and product
\[\qquad\qquad\quad C\mapsto C^*\qquad\textrm{and}\qquad(B,C)\mapsto(BC)^\precsim,\]
whenever $(B^*B)^\precsim=(CC^*)^\precsim$, in which case $(BC)^\precsim=(Bc)^\precsim$, for any $c\in C$.
\end{thm}

\begin{proof}
If $C$ is a coset then so is $C^*$ so the involution is well-defined on $\mathcal{C}(S)$.

Now take $B,C\in\mathcal{C}(S)$ with $(B^*B)^\precsim=(CC^*)^\precsim$ and $c\in C$.  We have $d\in C$ with $d\precsim c$ so, in particular, $d^*\prec cc^*$ and hence $dd^*\prec c^*c$.  Also
\[dd^*\in CC^*\subseteq C^\precsim C^{\precsim*}\subseteq(CC^*)^\precsim=(B^*B)^\precsim\]
so $dd^*\in\{cc^*\}^\succ\cap(B^*B)^\precsim$ and hence $(Bc)^\precsim$ is a coset, by \autoref{CaAtlas}.  For any other $a\in C$, we have $ac^*\in(CC^*)^\precsim=(B^*B)^\precsim$ so $ac^*B^*\subseteq B^{*\precsim}=B^*$, by \eqref{TUC<=TU}, and hence $(c^*B^*)^\precsim\subseteq(a^*B^*)^\precsim$, by \eqref{TUC=>TC}, i.e. $(Bc)^\precsim\subseteq(Ba)^\precsim$.  Switching $a$ and $c$ yields the reverse inclusion so $(Bc)^\precsim=(Ba)^\precsim$.  As $a$ was arbitrary, $(BC)^\precsim=(Bc)^\precsim$ is a coset so the product is also well-defined.

Again by \autoref{CaAtlas}, $Bc$ and $Bcc^*$ are atlases with $(Bcc^*)^\precsim=B^\precsim=B$ so
\[((BC)^\precsim(BC)^{\precsim*})^\precsim=((Bc)^\precsim(Bc)^{\precsim*})^\precsim=(Bcc^*B^*)^\precsim=((Bcc^*)^\precsim B^{*\precsim})^\precsim=(BB^*)^\precsim.\]
Thus the product of $A$ and $B$ is defined iff the product of $A$ and $(BC)^\precsim$ is defined.  Instead assuming $(AB)^\precsim$ and hence $(B^*A^*)^\precsim$ is valid product, this means the product of $C^*$ and $B^*$ is defined iff the product of $C^*$ and $(B^*A^*)^\precsim$ is defined, i.e. the product of $B$ and $C$ is defined iff the product of $(AB)^\precsim$ and $C$ is defined.  Then the products all equal $(ABC)^\precsim$, e.g.
\[(ABC)^\precsim=(A^\precsim B^\precsim C^\precsim)^\precsim\subseteq((AB)^\precsim C^\precsim)^\precsim\subseteq(ABC)^{\precsim\precsim}\subseteq(ABC)^\precsim.\]
This shows that the product is associative.

Finally, note that $(CC^*)^\precsim$ is a unit as, whenever $(BCC^*)^\precsim$ is a valid product, $(BCC^*)^\precsim=(Bcc^*)^\precsim=B^\precsim=B$, again by \eqref{Caa*}.
\end{proof}

In particular, $C\in\mathcal{C}(S)$ is a unit iff $C=(CC^*)^\precsim$. These unit cosets can be characterised in a number of ways.

\begin{prp}\label{IdempotentCosets}
For any non-empty $C\subseteq S$, the following are equivalent.
\begin{enumerate}
\item\label{UnitCoset} $C$ is a unit coset.
\item\label{CC*} $CC^*\subseteq C=C^\precsim$.
\item\label{*Subsemi} $C$ is a *-subsemigroup with $C=C^\precsim$.
\item\label{CosetS+} $C$ is a coset containing some $p\in |S|^2$.
\item\label{CosetE} $C$ is a coset containing some $e\in E$.
\end{enumerate}
\end{prp}

\begin{proof}\
\begin{itemize}
\item[\eqref{UnitCoset}$\Rightarrow$\eqref{CC*}] If $C$ is a unit coset then $CC^*\subseteq C^\precsim C^{\precsim*}\subseteq(CC^*)^\precsim=C$.
\item[\eqref{CC*}$\Rightarrow$\eqref{*Subsemi}] If $CC^*\subseteq C=C^\precsim$ then $C^*\subseteq C^{*\precsim}\subseteq(CC^*C^*)^\precsim\subseteq(CC^*)^\precsim\subseteq C^\precsim\subseteq C$.  Thus $C^*=C$ and $CC\subseteq CC^*\subseteq C$, i.e. $C$ is a *-subsemigroup.
\item[\eqref{*Subsemi}$\Rightarrow$\eqref{CosetS+}] If $CC\subseteq C=C^*=C^\precsim$ then $CC^*C\subseteq C$, so $C$ is a coset, and $|C|^2\subseteq C$.
\item[\eqref{CosetS+}$\Rightarrow$\eqref{CosetE}] Immediate from $|S|^2\subseteq E$.
\item[\eqref{CosetE}$\Rightarrow$\eqref{UnitCoset}] If $C$ is a coset and $e\in C\cap E\subseteq E^*$ then \eqref{EBelow} yields
\[C\subseteq C^\precsim\subseteq(Ce^*)^\precsim\subseteq(CC^*)^\precsim\subseteq(CC^*e)^\precsim\subseteq(CC^*C)^\precsim\subseteq C^\precsim\subseteq C.\]
Thus $C=(CC^*)^\precsim$ so $C$ is a unit coset.
\end{itemize}
\end{proof}

The canonical topology on $\mathcal{C}=\mathcal{C}(S)$ is generated by the sets $(\mathcal{C}_a)_{a\in S}$ given by
\[\mathcal{C}_a=\{C\in\mathcal{C}:a\in C\}.\]
In other words, we are taking $(\mathcal{C}_a)_{a\in S}$ as a subbasis for the topology on $\mathcal{C}$.  To get a basis, we consider the sets $\mathcal{C}_F$, for finite $F\subseteq S$, where
\[\mathcal{C}_F=\bigcap_{f\in F}\mathcal{C}_f=\{C\in\mathcal{C}:F\subseteq C\}.\]

\begin{thm}\label{EtaleCosets}
The coset groupoid $\mathcal{C}=\mathcal{C}(S)$ is \'etale.
\end{thm}

\begin{proof}
As $\mathcal{C}_a^*=\mathcal{C}_{a^*}$, the involution is certainly continuous.  Also, if $(BC)^\precsim\in\mathcal{C}_a$ then we must have $b\in B$ and $c\in C$ with $bc\precsim a$ so $(BC)^\precsim\in\mathcal{C}_b\mathcal{C}_c\subseteq\mathcal{C}_a$ and hence the product is also continuous.

It only remains to show that the range $r$ is an open map.  To see this, take $B\in\mathcal{C}_F$, for some finite $F\subseteq S$.  For each $f\in F$, we have $b_f\in B$ with $b_f\precsim f$.  Fix any $a\in F$ and let $b=b_a\precsim a$.  In particular, $b^*\prec aa^*$ and hence $bb^*\prec aa^*$.  Now any $C\in\bigcap_{f\in F}\mathcal{C}_{b_fb^*}\subseteq\mathcal{C}_{bb^*}$ is a unit, by \autoref{IdempotentCosets}, so it follows that $bb^*\in\{aa^*\}^\succ\cap C=\{aa^*\}^\succ\cap(C^*C)^\precsim$.  By \autoref{CaAtlas}, $Ca$ is an atlas with
\[r((Ca)^\precsim)=(Ca)^\precsim(Ca)^{\precsim*}=(Caa^*C^*)^\precsim=(CC^*)^\precsim=C.\]
For each $f\in F$, $b_fb^*\in C$ implies that $a_f\succsim b_fb^*a\in Ca$, as $b^*a\in E$, and hence $a_f\in(Ca)^\precsim$, i.e. $(Ca)^\precsim\in\mathcal{C}_F$.  As $C$ was arbitrary,
\[r(B)=(BB^*)^\precsim\in\bigcap_{f\in F}\mathcal{C}_{b_fb^*}\subseteq r[\mathcal{C}_F].\]
As $B$ was arbitrary, $r[\mathcal{C}_F]$ is open in $\mathcal{C}$.  As $F$ was arbitrary, $r$ is an open map.
\end{proof}

While $\mathcal{C}(S)$ is \'etale, it often has bad separation properties, e.g. if there is at least one non-empty coset properly containing another then $\mathcal{C}(S)$ is not even $T_1$, let alone Hausdorff.  One way of correcting this is to restrict to a subspace of $\mathcal{C}(S)$, e.g. consisting of special kinds of filters, which we investigate next.

\subsection{Filters}\label{Filters}

\begin{dfn}
We call non-empty $U\subseteq S$ a \emph{filter} if
\[\tag{Filter}\label{Filter}a,b\in U\qquad\Leftrightarrow\qquad\exists c\in U\ (c\precsim a,b).\]
\end{dfn}

In other words, a filter is a down-directed up-set, i.e. $U^\precsim\subseteq U$ and every finite $F\subseteq U$ has a lower $\precsim$-bound (including $F=\{a\}$, so $U\subseteq U^\precsim$, and $F=\emptyset$, so $U\neq\emptyset$).  One advantage of filters over cosets is that, once we know the elements below some point in the filter, we already know the filter.  Indeed, one can verify that
\[U\textrm{ is a filter}\qquad\Leftrightarrow\qquad\forall a\in U\ (U=(U\cap a^\succsim)^\precsim).\]

We will soon see that every filter is a coset, but there can certainly be non-filter cosets.  E.g. consider the situation in \autoref{Groups}, where $S$ is a group and $E$ is a normal subgroup.  Then the cosets are precisely the subsets of the form $aG$, for some $a\in S$ and subgroup $G$ containing $E$.  However, $aG$ here will be a filter (if and) only if $G=E$.  In particular, when $E=|S|^2=\{u\}(=$ the unit of $S)$, cosets come from arbitrary subgroups, while filters are precisely the singletons.  There is, however, a canonical unit filter associated to any non-empty coset $C$, namely $|C|^{2\precsim}$.

\begin{prp}\label{C2Filter}
If $C\subseteq S$ is a non-empty atlas then $|C|^2$ is down-directed and
\[C\textrm{ is down-directed}\quad\Leftrightarrow\quad C^*C\textrm{ is down-directed}\quad\Leftrightarrow\quad(C^*C)^\precsim=|C|^{2\precsim}.\]
\end{prp}

\begin{proof}
If $C$ is an atlas then $C^\precsim$ is a coset so, for any $a,b\in C\subseteq C^\precsim$, \autoref{*Directed} yields $c\in C$ with $c\precsim a$ and $c\prec b^*b$.  Thus $c^*c\precsim a^*a,b^*b$ so $|C|^2$ is down-directed.

If $C$ is down-directed then $C^*C$ is down-directed, by \eqref{ProductPreserving}.  Conversely, assume $C^*C$ is down-directed and take any $a,b\in C$.  As $C\subseteq C^\precsim$, we have $c,d\in C$ with $c\precsim a$ and $d\precsim b$.  As $C^*C$ is down-directed, we have $f\in C^*C$ with $f\precsim c^*c,c^*d$.  Note $f=c^*cc^*cf=a^*ac^*cc^*cf=a^*af$ so \eqref{TripleProduct} yields $af\precsim ac^*c,ac^*d$ and hence $af\precsim a,b$, as $ac^*\in E$.  But $af\in CC^*C\subseteq C^\precsim$, so we have $g\in C$ with $g\precsim af\precsim a,b$, showing that $C$ is down-directed.  This proves the first equivalence.

If $C$ is down-directed then, for any $a,b\in C$, we have $c\in C$ with $c\precsim a,b$ and hence $c^*c\precsim a^*b$, showing that $(C^*C)^\precsim\subseteq|C|^{2\precsim}\subseteq(C^*C)^\precsim$.  Conversely, if $(C^*C)^\precsim=|C|^{2\precsim}$ then, as we already known $|C|^2$ is down-directed, $C^*C$ is too.
\end{proof}

\begin{prp}\label{C2C*C2}
If $C\subseteq S$ is an atlas then $|C|^{2\prec}=|C^*C|^{2\prec}$.
\end{prp}

\begin{proof}
For any $b\in C$, we have $c\in C$ with $c\precsim b$ so $c\prec b^*b$ and hence $c^*cc^*c\prec b^*b$, showing that $|C|^2\subseteq|C^*C|^{2\prec}$.  On the other hand, for any $c,d\in C$, \autoref{*Directed} yields $b\in C$ with $b\precsim d$ and $b^*\prec cc^*$ so $bd^*\in E$ commutes with $cc^*\in|S|^2$ and hence $b(c^*d)^*c^*d=bd^*cc^*d=cc^*bd^*d=b$.  Thus $b\prec (c^*d)^*c^*d$ and hence $b^*b\prec (c^*d)^*c^*d$ so $|C^*C|^2\subseteq|C|^{2\prec}$.
\end{proof}

\begin{prp}\label{FilterSubgroupoid}
The filters form an ideal in the coset groupoid $\mathcal{C}(S)$.
\end{prp}

\begin{proof}
Assume $U\subseteq S$ is a filter.  For any $a,b,c\in U$, take $u,v\in U$ with $u\precsim v\precsim a,b,c$.  By \eqref{TripleProduct} and \eqref{ProductPreserving}, $u\precsim vv^*v\precsim ab^*c$
so $ab^*c\in U^{\precsim\precsim}=U$.  Thus $UU^*U\subseteq U$ so $U$ is a coset, i.e. all filters are cosets.  Also $U$ is a filter iff $U^*$ is a filter, by \eqref{*Invariance}.  And a coset $U$ is a filter iff $(U^*U)^\precsim$ is a filter, by the first equivalence in \autoref{C2Filter}, so the filters do indeed form an ideal.
\end{proof}

The second equivalence in \autoref{C2Filter} shows that unit filters are generated by their *-squares.  In fact, for a unit coset to be a filter, it suffices that it is generated by the elements in $E$.  In other words, while cosets containing a single element of $E$ are necessarily units, if they contain sufficiently many elements of $E$ then they must also be filters.

\begin{prp}
The unit filters are the non-empty $E$-generated cosets, i.e.
\[U\textrm{ is a unit filter}\qquad\Leftrightarrow\qquad U=(U\cap E)^\precsim\textrm{ is a coset}.\]
\end{prp}

\begin{proof}
If $U$ is a unit filter then $U$ contains some $e\in E$.  Any other $a\in U$ shares a lower bound in $U$ with $e$, which must also be in $E$, by \eqref{EDownClosed}, i.e. $U=(U\cap E)^\precsim$.

Conversely, if $U$ is a unit coset, we claim $U\cap E$ is down-directed.  To see this, take $e,f\in U\cap E$, so we have $e',f'\in U$ with $e'\precsim e$ and $f'\precsim f$ and hence $e',f'\in E$, by \eqref{EDownClosed}.  By \autoref{IdempotentCosets}, $U$ is a subsemigroup so $e'f'\in U\cap E$ and $e'f'\precsim e,f$, by \eqref{EBelow}.  This proves the claim, from which it follows $U=(U\cap E)^\precsim$ is a filter.
\end{proof}

\begin{dfn}
A maximal proper filter is an \emph{ultrafilter}.
\end{dfn}

In other words, a filter $U\subsetneqq S$ is an ultrafilter if $U$ is not contained in any larger filter, except possibly $S$ itself.

\begin{prp}\label{UltrafilterIdeal}
The ultrafilters form an ideal in the filter groupoid.
\end{prp}

\begin{proof}
By \eqref{*Invariance}, $U$ is an ultrafilter iff $U^*$ is an ultrafilter.  Thus it suffices to show that $(UV)^\precsim$ is an ultrafilter whenever $U$ is a filter and $V$ is an ultrafilter with $(U^*U)^\precsim=(VV^*)^\precsim$.  Then $U=U^\precsim\neq\emptyset$ and $U^*UV\subseteq V\neq S$ so \eqref{TUC=>Proper} yields $(UV)^\precsim\neq S$.  So if $(UV)^\precsim$ were not an ultrafilter, we would have a proper filter $W\supsetneqq(UV)^\precsim$.  Then we could take $w\in W\setminus(UV)^\precsim$ and $x\in W$ with $x\precsim w$.  For any $u\in U$, note that that $u^*x\notin V$ \textendash\, otherwise we would have $w\succsim uu^*x\in UV$ so $w\in(UV)^\precsim$, contradicting our choice of $w$.  But this means that $V\subsetneqq(U^*W)^\precsim$, as $V=(U^*UV)^\precsim\subseteq(U^*W)^\precsim$.  By \eqref{ProductPreserving}, $(U^*W)^\precsim$ is a filter, as both $U$ and $W$ are down-directed.  Moreover, taking $T=UVV^*$, we have $T^\precsim=U\neq\emptyset$ and $TU^*W\subseteq UVV^*U^*W\subseteq WW^*W\subseteq W$ so $(U^*W)^\precsim\neq S$, again by \eqref{TUC=>Proper}.  But this contradicts the maximality of $V$, so we are done.
\end{proof}

\begin{dfn}
The \emph{Weyl groupoid} is the set of ultrafilters
\[\mathcal{U}(S)=\{U\subseteq S:U\textrm{ is an ultrafilter}\}\]
considered as a topological subgroupoid of the coset groupoid $\mathcal{C}(S)$.
\end{dfn}

So the topology on $\mathcal{U}(S)$ is generated by the sets
\[\mathcal{U}_a=\{V\in\mathcal{U}:a\in V\}.\]
As every $V\in\mathcal{U}(S)$ is a filter, $(\mathcal{U}_a)_{a\in S}$ is already a basis for this topology.

\begin{rmk}\label{TightWeyl}
By \cite[Theorem 5.15]{LawsonLenz2013}, when $S$ is an inverse semigroup and $E=|S|^2$, the Weyl groupoid coincides with the Lawson-Lenz version of Exel's tight groupoid if (and only if) $S$ satisfies a certain `trapping' condition similar to \autoref{TrappingCorollary} below.  Indeed, while it would be formally possible to define tight filters along the lines of \cite{BiceStarling2020HTight}, these would turn out to be the same as ultrafilters in the *-subsemigroups of *-rings that we are primarily concerned with here.

On the other hand, if $A$ is a C*-algebra with Cartan subalgebra $C$ and we take $E=C^1(=$ the unit ball of $C)$ and $S=C_*^1(=$ the contractive normalisers of $C)$ then $\mathcal{U}(S)$ coincides with the Kumjian-Renault Weyl groupoid.  Indeed, the Kumjian-Renault construction yields a groupoid $G$ on which $S$ gets represented as continuous functions supported on bisections, and hence $G$ must coincide with $\mathcal{U}(S)$, by \autoref{GroupoidRecovery}.
\end{rmk}

Up till now, the Weyl groupoid has only been considered in the context of C*-algebras.  In fact, several of its important properties remain valid in more general *-rings, and some can even be proved at the *-semigroup level, like the following.

\begin{thm}\label{EtaleGroupoid}
The Weyl groupoid $\mathcal{U}(S)$ is \'etale.
\end{thm}

\begin{proof}
This is immediate from \autoref{OpenIdeals}, \autoref{EtaleCosets}, \autoref{FilterSubgroupoid} and \autoref{UltrafilterIdeal}.
\end{proof}

\begin{rmk}
We could have restricted our attention to ultrafilters from the start, as these constitute the Weyl groupoid that we are primarily interested in.  However, as shown in \autoref{EtaleCosets}, it takes little extra effort to show that the larger coset groupoid is \'etale, from which \autoref{EtaleGroupoid} is immediate.  We suspect cosets will also play a more important role in future work examining Weyl-like groupoids with better separation and functorial properties, just as they did with tight groupoids in \cite{BiceStarling2020HTight}.
\end{rmk}

Recall that $0\in S$ is \emph{absorbing} $0a=0=a0$, for all $a\in S$.  Equivalently, $0$ is absorbing iff $0^\prec=S$.  Indeed $0^\prec=S$ means that $0=0a$, for all $a\in S$, but then $0=00^*=(00^*)^*=0^*$ and hence $0a=(a^*0^*)^*=(a^*0)^*=0^*=0$, for all $a\in S$.  This also shows that $0=00^*\in |S|^2\subseteq E$ and hence $S=0^\sim=0^\precsim$.

\begin{prp}\label{HausdorffUnitSpace}
If $0\in S$ is absorbing, the unit space $\mathcal{U}^0=\mathcal{U}(S)^0$ is Hausdorff.
\end{prp}

\begin{proof}
Take any distinct $U,V\in\mathcal{U}^0$.  By \eqref{ProductPreserving}, $UV$ is down-directed so $(UV)^\precsim$ is a filter.  As $V\in\mathcal{U}^0$, we have $e\in V\cap E$, by \autoref{IdempotentCosets}, and hence $U\subseteq(Ue)^\precsim\subseteq(UV)^\precsim$.  Likewise, $V\subseteq(UV)^\precsim$.  If we had $U=(UV)^\precsim$ then we would have $U\supseteq V$ and hence $U=V$, as $V$ is a ultrafilter, contradicting our choice of $U$ and $V$.  Thus $U\neq(UV)^\precsim$ and hence $(UV)^\precsim=S$, as $U$ is an ultrafilter.  In particular, $0\in(UV)^\precsim$, so we have $u\in U$ and $v\in V=V^*$ with $uv\precsim0$ and hence $uv=uv00=0$.  But this means $\mathcal{U}_u\cap\mathcal{U}_{v^*}=\emptyset$, as any $W\in\mathcal{U}_u\cap\mathcal{U}_{v^*}$ would contain $w\precsim u,v^*$ and hence $w=wvv^*=wu^*uvv^*=0$, which implies $W=W^\precsim\supseteq0^\precsim=S$, contradicting the fact $W$ is an ultrafilter.  Thus $\mathcal{U}^0$ is Hausdorff.
\end{proof}

Even when $S$ is an inverse semigroup, it is known that entire Weyl/tight groupoid $\mathcal{U}(S)$ may fail to be Hausdorff.  But by \autoref{EtaleGroupoid}, $\mathcal{U}(S)$ is \'etale and hence locally homeomorphic to $\mathcal{U}^0$ so, by \autoref{HausdorffUnitSpace}, $\mathcal{U}(S)$ is at least locally Hausdorff.

\subsection{Actions}\label{Actions}

Recall that an \emph{action} of $S$ on a set $X$ is a map $(s,x)\mapsto sx$ from $S\times X$ to $X$ satisfying the natural associativity axiom w.r.t. the product in $S$, i.e.
\[(st)x=s(tx).\]
Actions on uniform spaces lead to a natural bundle construction over the Weyl groupoid $\mathcal{U}(S)$, generalising the Weyl line bundle considered in the context of C*-algebras (see \cite[Definition 16.13]{ExelPitts2019}).

First recall that the \emph{composition} $\mathsf{P}\circ\mathsf{Q}$ of relations $\mathsf{P},\mathsf{Q}\subseteq X\times X$ is defined by
\[\mathsf{P}\circ\mathsf{Q}=\{(x,y):\exists z\in X\ ((x,z)\in\mathsf{P}\textrm{ and }(z,y)\in\mathsf{Q})\}.\]
In the usual infix notation for relations, this can equivalently be written as
\[x\mathrel{(\mathsf{P}\circ\mathsf{Q})}y\qquad\Leftrightarrow\qquad\exists z\in X\ (x\mathrel{\mathsf{P}}z\mathrel{\mathsf{Q}}y).\]
A \emph{uniformity} on $X$ is a $\subseteq$-filter $\mathfrak{F}\subseteq\mathcal{P}(X\times X)$ of symmetric reflexive relations with
\[\tag{Divisible}\mathsf{R}\in\mathfrak{F}\qquad\Rightarrow\qquad\exists\mathsf{Q}\in\mathfrak{F}\ (\mathsf{Q}\circ\mathsf{Q}\subseteq\mathsf{R}).\]
A \emph{uniform space} is simply a set $X$ together with a uniformity $\mathfrak{F}$ on $X$.

Considering $S$ as a discrete uniform space, the canonical product uniformity on $S\times X$ is generated by the relations $\mathsf{R}'=\{((s,x),(s,y)):(x,y)\in\mathsf{R}\}$, for $\mathsf{R}\in\mathfrak{F}$.  We say the action of $S$ on $X$ is uniformly continuous if $(s,x)\mapsto sx$ is uniformly continuous with respect to this product uniformity.  More explicitly this means that, with $S\mathsf{Q}=\{(sx,sx):(x,y)\in\mathsf{Q}\}$,
\begin{equation}\label{UniformAction}
\mathsf{R}\in\mathfrak{F}\qquad\Rightarrow\qquad\exists\mathsf{Q}\in\mathfrak{F}\ (S\mathsf{Q}\subseteq\mathsf{R}).
\end{equation}

To avoid unnecessary repetition, we make the following assumption throughout the rest of this section (it will no longer be needed in the following section).

\begin{ass}\label{UniformityAssumption}
\textbf{We have a uniformly continuous action of $S$ on a uniform space $X$ and an $E$-equivariant map $\Psi$ on $X$, i.e. for all $e\in E$,}
\[\label{EEquivariant}\tag{$E$-Equivariant}\Psi(ea)=e\Psi(a).\]
\end{ass}

Given an arbitrary action of $S$ on a set $X$, we could just consider the discrete uniformity on $X$ and take $\Psi$ to be the identity map.  In this case, our Weyl bundle construction would be reminiscent of the germs used in Exel's original version of the tight groupoid \textendash\, see \cite[Definition 4.6]{Exel2008}.

However, the motivating situation we have in mind is that $X$ is a C*-algebra or, more generally, a Hilbert $A$-module, for some C*-algebra $A$.  Moreover, we imagine $S$ is a *-subsemigroup of the unit ball of $A$ so the module structure immediately yields an action of $S$ on $X$.   We further imagine $\Psi$ is a conditional expectation onto a Cartan subalgebra or submodule and the uniformity on $X$ is generated by the relations $(\equiv^\delta)_{\delta>0}$ where
\begin{equation}\label{deltaEquiv}
x\equiv^\delta y\qquad\Leftrightarrow\qquad\|x-y\|<\delta.
\end{equation}
As $\|sx-sy\|\leq\|s\|\|x-y\|=\|x-y\|$, for any $s$ in the unit ball of $A$, the uniform continuity condition in \autoref{UniformityAssumption} is indeed satisfied, even with $\mathsf{Q}=\mathsf{R}$ in \eqref{UniformAction}.

\begin{dfn}
For every $U\subseteq S$ and $\mathsf{R}\in\mathfrak{F}$, we define a relation $\mathsf{R}_U$ on $X$ by
\begin{align*}
x\mathrel{\mathsf{R}_U}y\qquad&\Leftrightarrow\qquad\exists u\in U\ (\Psi(u^*x)\mathrel{\mathsf{R}}\Psi(u^*y)).\\
\intertext{We then define a relation $\equiv_U$ on $A$ by}
x\equiv_Uy\qquad&\Leftrightarrow\qquad\forall\mathsf{R}\in\mathfrak{F}\ (x\mathrel{\mathsf{R}_U}y).
\end{align*}
\end{dfn}

\begin{prp}
If $U\subseteq S$ is down-directed then $\equiv_U$ is an equivalence relation.
\end{prp}

\begin{proof}
Note $\equiv_U$ is symmetric and reflexive, as each $\mathsf{R}\in\mathfrak{F}$ is.  For transitivity, assume $x\equiv_Uy\equiv_Uz$.  As $\mathfrak{F}$ is a uniformity, for any $\mathsf{R}\in\mathfrak{F}$, we have $\mathsf{Q}\in\mathfrak{F}$ with $\mathsf{Q}\circ\mathsf{Q}\subseteq\mathsf{R}$.  As the action of $S$ on $X$ is uniformly continuous, we also have $\mathsf{R}\in\mathfrak{F}$ with $S\mathsf{R}\subseteq\mathsf{Q}$.  As $x\equiv_Uy\equiv_Uz$, we have $s,t\in U$ with $\Psi(s^*x)\mathrel{\mathsf{R}}\Psi(s^*y)$ and $\Psi(t^*y)\mathrel{\mathsf{R}}\Psi(t^*z)$.  As $U$ is down-directed, we have $u\in U$ with $u\precsim s,t$.  As $\Psi(ed)=e\Psi(d)$, for all $e\in E$,
\[\Psi(u^*x)=u^*s\Psi(s^*x)\mathrel{\mathsf{Q}}u^*s\Psi(s^*y)=\Psi(u^*y)=u^*t\Psi(t^*y)\mathrel{\mathsf{Q}}u^*t\Psi(t^*z)=\Psi(u^*z)\]
and hence $\Psi(u^*x)\mathrel{\mathsf{R}}\Psi(u^*z)$.  This shows that $\equiv_U$ is indeed transitive.
\end{proof}

\begin{dfn}\label{WeylBundle}
The \emph{Weyl bundle} is given by
\[\mathcal{W}=\mathcal{W}(X)=\{(U,x^{\equiv_U}):U\in\mathcal{U}(S)\textrm{ and }x\in X\}.\]
\end{dfn}

In other words, the Weyl bundle consists of ordered pairs $(U,Y)$ where $U$ is an ultrafilter in $S$ and $Y$ is a $\equiv_U$-equivalence class in $X$.  In the usual context considered in C*-algebras, this provides an alternative construction of the Weyl line bundle (see \cite[Definition 16.13]{ExelPitts2019}), at least when we give it the appropriate topological and algebraic structure.  In this section we focus on the topology.

First, for $x\in X$, $s\in S$ and $\mathsf{R}\in\mathfrak{F}$, define $x_s^\mathsf{R}\subseteq\mathcal{W}$ by
\[x_s^\mathsf{R}=\{(U,y^{\equiv_U})\in\mathcal{W}:s\in U\textrm{ and }x\mathrel{\mathsf{R}_U}y\}.\]
Note $(U,x^{\equiv_U})\in x_s^\mathsf{R}$, for all $(U,x^{\equiv_U})\in\mathcal{W}$ and $s\in U$, as each $\mathsf{R}\in\mathfrak{F}$ is reflexive.
Take
\[\mathcal{N}_{(U,Y)}=\{y_s^\mathsf{R}:s\in U,y\in Y\textrm{ and }\mathsf{R}\in\mathfrak{F}\}\]
as a neighbourhood subbase at $(U,Y)$.  Actually, it turns out this defines a neighbourhood base at $(U,Y)$, even when we fix some $y\in Y$.

\begin{prp}\label{SingleNBase}
Any $(U,x^{\equiv_U})\in\mathcal{W}$ has a neighbourhood base of the form
\[\mathcal{N}=\{x_s^\mathsf{R}:s\in U\textrm{ and }\mathsf{R}\in\mathfrak{F}\}.\]
\end{prp}

\begin{proof}
Assume we have a finite collection $x_{1s_1}^{\mathsf{R}_1},\cdots,x_{ks_k}^{\mathsf{R}_k}$ of subbasic neighbourhoods of $(U,x^{\equiv_U})$.  Take $\mathsf{Q}\in\mathfrak{F}$ with $\mathsf{Q}\circ\mathsf{Q}\subseteq\mathsf{R}_1\cap\cdots\cap\mathsf{R}_k$ and take $\mathsf{R}\in\mathfrak{F}$ with $S\mathsf{R}\subseteq\mathsf{Q}$.  For all $j\leq k$, $x\equiv_Ux_j$ so we have $u_j\in U$ with $\Psi(u_j^*x)\mathrel{\mathsf{R}}\Psi(u_j^*x_j)$.  Taking $u\in U$ with $u\precsim u_1,\cdots,u_k,s_1,\cdots,s_k$, we claim that $x_u^\mathsf{R}\subseteq x_{1s_1}^{\mathsf{R}_1}\cap\cdots\cap x_{ks_k}^{\mathsf{R}_k}$.  To see this take $(V,Y)\in x_u^\mathsf{R}$, so $u\in V$ and $x\mathrel{\mathsf{R}_V}y$, for some $y\in Y$, which means we have $v\in V$ with $\Psi(v^*x)\mathrel{\mathsf{R}}\Psi(v^*y)$.  Taking $w\in V$ with $w\precsim u,v$, we see that
\[\Psi(w^*x_j)=w^*u_j\Psi(u_j^*x_j)\mathrel{\mathsf{Q}}w^*u_j\Psi(u_j^*x)=\Psi(w^*x)=w^*v\Psi(v^*x)\mathrel{\mathsf{Q}}w^*v\Psi(v^*y)=\Psi(w^*y)\]
and hence $\Psi(w^*x_j)\mathrel{\mathsf{R}_j}\Psi(w^*y)$.  As $w\in V$, this means that $x_j\mathrel{\mathsf{R}_{jV}}y$.  As $y\in Y$ and $s_j\succsim u\in V$, this means that $(V,Y)\in x_{js_j}^{\mathsf{R}_j}$, as required.
\end{proof}

So $N\subseteq\mathcal{W}$ is a neighbourhood of $(U,x^{\equiv_U})$ iff $x_u^\mathsf{R}\subseteq N$, for some $u\in U$ and $\mathsf{R}\in\mathfrak{F}$.  As usual, we define the interior of any $N\subseteq\mathcal{W}$ to be the points having $N$ as a neighbourhood.  To show that the neighbourhoods define a topology, rather than just a pretopology, we have to show that the interior of any neighbourhood is again a neighbourhood of the same point.

\begin{prp}
The neighbourhoods above define a topology on $\mathcal{W}$.
\end{prp}

\begin{proof}
Take a basic neighbourhood $x_u^\mathsf{R}$ of $(U,x^{\equiv_U})\in\mathcal{W}$, so $u\in U$ and $\mathsf{R}\in\mathfrak{F}$.  Take $\mathsf{P},\mathsf{Q}\in\mathfrak{F}$ with $\mathsf{Q}\circ\mathsf{Q}\subseteq\mathsf{R}$ and $S\mathsf{P}\subseteq\mathsf{Q}$.  We claim that $x_u^\mathsf{P}$ is in the interior of $x_u^\mathsf{R}$.  To see this, take $(V,Y)\in x_u^\mathsf{P}$, so we have $y\in Y$ with $x\mathrel{\mathsf{P}_V}y$, which means we have $v\in V$ with $\Psi(v^*x)\mathrel{\mathsf{P}}\Psi(v^*y)$.  Take $t\in V$ with $t\precsim u,v$.  To show $(V,Y)$ is in the interior of $x_u^\mathsf{P}$, it suffices to show that $y^\mathsf{P}_t\subseteq x^\mathsf{R}_u$.  If $(W,Z)\in y^\mathsf{P}_t$ then we have $z\in Z$ with $y\mathrel{\mathsf{P}_W}z$, which means we have $w\in W$ with $\Psi(w^*y)\mathrel{\mathsf{P}}\Psi(w^*z)$.  Taking $s\in W$ with $s\precsim t,w$, we get
\[\Psi(s^*x)=s^*v\Psi(v^*x)\mathrel{\mathsf{Q}}s^*v\Psi(v^*y)=\Psi(s^*y)=s^*w\Psi(w^*y)\mathrel{\mathsf{Q}}s^*w\Psi(w^*z)=\Psi(s^*z).\]
Thus $\Psi(s^*x)\mathrel{\mathsf{R}}\Psi(s^*z)$ and hence $x\mathrel{\mathsf{R}_W}z$ so $(W,Z)\in x^\mathsf{R}_u$.  As $(W,Z)$ was arbitrary, $y^\mathsf{P}_t\subseteq x^\mathsf{R}_u$.  As $(V,Y)$ was arbitrary, this shows that $x_u^\mathsf{P}$ is in the interior of $x_u^\mathsf{R}$ and hence this interior is still a neighbourhood of $(U,x^{\equiv_U})$, as required.
\end{proof}

Next we note that the Weyl bundle is indeed a bundle over the Weyl groupoid.

\begin{prp}\label{GroupoidBundle}
The projection $\pi((U,Y))=U$ is an open continuous map from the Weyl bundle $\mathcal{W}=\mathcal{W}(X)$ onto the Weyl groupoid $\mathcal{U}=\mathcal{U}(S)$.
\end{prp}

\begin{proof}
To see that $\pi$ is continuous just note $\pi^{-1}[\mathcal{U}_s]$ is open, for all $s\in S$, as
\[\pi^{-1}[\mathcal{U}_s]=\bigcup_{x\in X}x^\mathsf{R}_s,\]
for any $\mathsf{R}\in\mathfrak{F}$.  On the other hand, for any $x\in X$, $s\in S$ and $\mathsf{R}\in\mathfrak{F}$,
\[\pi[x_s^\mathsf{R}]=\mathcal{U}_s.\]
Indeed, for any $U\in\mathcal{U}_s$, $(U,x^{\equiv_U})\in x_s^\mathsf{R}$ and $\pi((U,x^{\equiv_U}))=U$.  Thus $\pi$ is an open map onto $\mathcal{U}$.
\end{proof}

\begin{prp}\label{Representation}
Every $x\in X$ defines a continuous section $\widehat{x}$ of $\mathcal{W}$ given by
\[\widehat{x}(U)=(U,x^{\equiv_U}).\]
\end{prp}

\begin{proof}
Take a basic neighbourhood $x^\mathsf{R}_s$ of $\widehat{x}(U)=(U,x^{\equiv_U})$ in $\mathcal{W}$.  Then $\mathcal{U}_s$ is a neighbourhood of $U$ in $\mathcal{U}$ and, for any $V\in\mathcal{U}_s$, we see that $\widehat{x}(V)=(V,x^{\equiv_V})\in x^\mathsf{R}_s$, i.e. $\mathcal{U}_s\subseteq\widehat{x}^{-1}[x^\mathsf{R}_s]$.  Thus $\widehat{x}$ is continuous.
\end{proof}

Lastly, we note $\mathcal{W}(X)$ may not be Hausdorff, even when the base space $\mathcal{U}(S)$ is.

\begin{xpl}\label{NonHausdorffBundle}
Consider the $\mathbb{Z}$-sequences of elements in the complex unit disk whose left and right subsequences converge to the same element, i.e. let
\[S=E=\{(\alpha_n)_{n\in\mathbb{Z}}\subseteq\mathbb{D}:\exists\alpha\in\mathbb{D}\ (\lim_{n\rightarrow\infty}\alpha_n=\alpha=\lim_{n\rightarrow\infty}\alpha_{-n})\}.\]
Then $S$ acts via pointwise multiplication on complex $\mathbb{Z}$-sequences whose left and right subsequences converge but this time to possibly different elements, i.e. we let
\[X=\{(\alpha_n)_{n\in\mathbb{Z}}\subseteq\mathbb{C}:\exists\beta,\gamma\in\mathbb{C}\ (\lim_{n\rightarrow\infty}\alpha_n=\beta\textrm{ and }\lim_{n\rightarrow\infty}\alpha_{-n}=\gamma)\}.\]
The Weyl groupoid of $S$ is then just a space, namely the one-point compactification of $\mathbb{Z}$.  More precisely, $\mathcal{U}(S)$ consists of ultrafilters $U_k=\{(\alpha_n)\in S:|\alpha_k|=1\}$, for $k\in\mathbb{Z}$, together with an extra point at infinity
\[U_\infty=\{(\alpha_n)\in S:\exists\textrm{ finite }F\subseteq\mathbb{Z}\ \forall k\in\mathbb{Z}\setminus F\ (|\alpha_k|=1)\}.\]
The Weyl bundle of $X$ has fibres we can identify with $\mathbb{C}$ at each $U_k$, but to $\mathbb{C}\times\mathbb{C}$ at $U_\infty$.  More precisely, $\mathcal{W}_{U_k}=\{(U_k,Y_k^\alpha):\alpha\in\mathbb{C}\}$ where
\[Y^\alpha_k=\{(\alpha_n)_{n\in\mathbb{Z}}\in X:\alpha_k=\alpha\},\]
for $k\in\mathbb{Z}$, but $\mathcal{W}_{U_\infty}=\{(U_\infty,Z^\beta_\gamma):\beta,\gamma\in\mathbb{C}\}$ where
\[Z^\beta_\gamma=\{(\alpha_n)_{n\in\mathbb{Z}}\in X:\lim_{n\rightarrow\infty}\alpha_n=\beta\textrm{ and }\lim_{n\rightarrow\infty}\alpha_{-n}=\gamma\}.\]
Let $x^\beta_\gamma\in S$ be the $\mathbb{Z}$-sequence taking the constant values $\beta$ and $\gamma$ on positive and negative integers respectively and let $s_F\in S$, for $F\subseteq\mathbb{Z}$, be the $\mathbb{Z}$-sequence which is $0$ on $F$ and $1$ elsewhere.  Then $(x^\beta_\gamma)^{\equiv^\delta}_{s_F}$, for finite $F\subseteq\mathbb{Z}$ and $\delta>0$, is a neighbourhood base of $Z^\beta_\gamma$ (where $(x_n)_{n\in\mathbb{Z}}\equiv^\delta(y_n)_{n\in\mathbb{Z}}$ means $\sup_{n\in\mathbb{Z}}|x_n-y_n|<\delta$).  As $F$ is finite, $(x^\beta_\gamma)^{\equiv^\delta}_{s_F}$ always contains $Y^\beta_k$, for some positive $k\in\mathbb{Z}$.  Thus $Y^\beta_n\rightarrow Z^\beta_\gamma$, as $n\rightarrow\infty$, for all $\gamma\in\mathbb{C}$.  In particular, limits are not unique so $\mathcal{W}(X)$ is not Hausdorff.
\end{xpl}

\section{*-Rings}\label{*Rings}

Further properties of the Weyl groupoid $\mathcal{U}(S)$ require our our Weyl *-semigroup $S$ to be sitting inside an ambient \emph{*-ring} \textendash\, a structure of the form $(A,+,\cdot,{}^*,0)$ where $A$ is a ring that is simultaneously a *-semigroup with respect to both sums and products.  Our primary goal is to show that, under suitable assumptions, $\mathcal{U}(S)$ is locally compact \textendash\, see \autoref{WeylLocallyCompact} below \textendash\, but to do this we first need to develop some general *-ring theory.

\begin{ass}
\textbf{We are given a *-ring $A$.}
\end{ass}

Note we do not assume that $A$ is necessarily unital.

\subsection{Norms}\label{Norms}

The key motivating examples of *-rings are of course C*-algebras.  In C*-algebras, it is well known that the algebra structure alone determines the norm via the spectrum.  Even without the scalars, we can define norm-like functions on any *-ring by instead using the order structure obtained from the *-squares.

Denote the additive subsemigroup generated by the *-squares of any $B\subseteq A$ by
\begin{align*}
B_\Sigma&=\{p_1+\cdots+p_k:p_1,\cdots,p_k\in|B|^2\}.\\
\intertext{Using $A_\Sigma$, we define functions from $A$ to $[0,\infty]$ (taking $\inf\emptyset=\infty$) by}
\lceil a\rceil&=\sup_{b\in A}\inf\{m/n:mbb^*-nbab^*\in A_\Sigma\}.\\
\|a\|&=\sqrt{\lceil aa^*\rceil\vee\lceil a^*a\rceil}.
\end{align*}
As $bA_\Sigma b^*\subseteq A_\Sigma$, when $A$ is unital, $\lceil\cdot\rceil$ reduces to
\[\lceil a\rceil=\inf\{m/n:m1-na\in A_\Sigma\}.\]
When $A$ is a $\mathbb{Q}$-algebra, this can be rewritten as $\lceil a\rceil=\inf\{r\in\mathbb{Q}:r1-a\in A_\Sigma\}$.  And when $A$ is a C*-algebra, $\lceil a\rceil=\|a_+\|$ (when $a\in A_\mathrm{sa}$, otherwise $\lceil a\rceil=\infty$) so $A_+=|A|^2=A_\Sigma=\{a\in A:\lceil-a\rceil=0\}$ and $\|a\|=\sqrt{\lceil aa^*\rceil}$ is the usual norm on $a$.

In a general *-ring, $\lceil\cdot\rceil$ is still an asymmetric seminorm (`asymmetric' here refers to the fact that we can have $\lceil-a\rceil\neq\lceil a\rceil$).

\begin{prp}
For all $a,b\in A$ and $n\in\mathbb{Z}_+$,
\begin{align}
\label{+Homogeneous}\tag{$+$-Homogeneous}\lceil na\rceil&=n\lceil a\rceil.\\
\label{Subadditive}\tag{Subadditive}\lceil a+b\rceil&\leq\lceil a\rceil+\lceil b\rceil.\\
\label{*Submultiplicative}\tag{*-Submultiplicative}\lceil bab^*\rceil&\leq\lceil bb^*\rceil\lceil a\rceil.
\end{align}
\end{prp}

\begin{proof}\
\begin{itemize}
\item[\eqref{+Homogeneous}] For $n=0$, note $bb^*-kb0b^*=bb^*\in|A|^2\subseteq A_\Sigma$, for all $k\in\mathbb{N}$ and $b\in A$, so
\[\lceil 0a\rceil=\lceil0\rceil=\inf_{k\in\mathbb{N}}1/k=0=0\lceil a\rceil.\]
For $n\in\mathbb{N}$, note that, for any $b\in A$ and $\epsilon>0$, we can find $j,k\in\mathbb{N}$ with $j/k\leq\lceil a\rceil+\epsilon$ and $jbb^*-kbab^*\in A_\Sigma$ and hence
\[jnbb^*-kb(na)b^*\in nA_\Sigma\subseteq A_\Sigma.\]
As $jn/k\leq n(\lceil a\rceil+\epsilon)$ and $b$ and $\epsilon$ were arbitrary, $\lceil na\rceil\leq n\lceil a\rceil$.  Conversely, we can find $j,k\in\mathbb{N}$ with $j/k\leq\lceil na\rceil+\epsilon$ and $jbb^*-knbab^*\in A_\Sigma$.  As $j/(kn)\leq\frac{1}{n}(\lceil na\rceil+\epsilon)$ we likewise get $\lceil a\rceil\leq\frac{1}{n}[na]$, i.e. $n\lceil a\rceil\leq\lceil na\rceil$.

\item[\eqref{Subadditive}] Take any $a,b\in A$.  For any $\epsilon>0$ and $c\in A$, we have $j,k\in\mathbb{N}$ with $j/k\leq\lceil a\rceil+\epsilon$ and $jcc^*-kcac^*\in A_\Sigma$.  Likewise, we have $m,n\in\mathbb{N}$ with $m/n\leq\lceil b\rceil+\epsilon$ and $mcc^*-ncbc^*\in A_\Sigma$.  Thus
\[(jn+km)cc^*-knc(a+b)c^*=n(jcc^*-kcac^*)+k(mcc^*-ncbc^*)\in nA_\Sigma+kA_\Sigma\subseteq A_\Sigma.\]
As $(jn+km)/kn=j/k+m/n\leq\lceil a\rceil+\lceil b\rceil+2\epsilon$ and $c$ and $\epsilon$ were arbitrary, it follows that $\lceil a+b\rceil\leq\lceil a\rceil+\lceil b\rceil$.

\item[\eqref{*Submultiplicative}] For any $\epsilon>0$, we have $j,k\in\mathbb{N}$ with $j/k\leq\lceil a\rceil+\epsilon$ and $jbb^*-kbab^*\in A_\Sigma$.  Likewise, for any $c\in A$, we have $m,n\in\mathbb{N}$ with $m/n\leq\lceil bb^*\rceil+\epsilon$ and $mcc^*-ncbb^*c^*\in A_\Sigma$ so
\[jmcc^*-kncbab^*c^*=jmcc^*-jncbb^*c^*+jncbb^*c^*-kncbab^*c^*\in jA_\Sigma+ncA_\Sigma c^*\subseteq A_\Sigma.\]
As $jm/(kn)=(j/k)(m/n)\leq(\lceil a\rceil+\epsilon)(\lceil bb^*\rceil+\epsilon)$, we are done.
\end{itemize}
\end{proof}

It follows that $\|\cdot\|$ is a *-invariant quasiseminorm (`quasi-' here refers to the fact that subadditivity only holds modulo some fixed factor, in this case $\sqrt{2}$).

\begin{cor}
For all $a,b\in A$ and $n\in\mathbb{Z}$,
\begin{align}
\label{*Invariant}\tag{*-Invariant}\|a^*\|&=\|a\|.\\
\label{Homogeneous}\tag{Homogeneous}\|na\|&=|n|\|a\|.\\
\label{Submultiplicative}\tag{Submultiplicative}\|ab\|&\leq\|a\|\|b\|.\\
\label{2Subadditive}\tag{$\sqrt{2}$-Subadditive}\tfrac{1}{\sqrt{2}}\|a+b\|&\leq\|a\|+\|b\|.
\end{align}
\end{cor}

\begin{proof}\
\begin{itemize}
\item[\eqref{*Invariant}] Immediate from $a^*a^{**}=a^*a$ and $a^{**}a^*=a^*a$.

\item[\eqref{Homogeneous}] Just note that \eqref{+Homogeneous} yields
\[\|na\|=\sqrt{\lceil n^2aa^*\rceil\vee\lceil n^2a^*a\rceil}=\sqrt{n^2(\lceil aa^*\rceil\vee\lceil a^*a\rceil)}=|n|\|a\|.\]

\item[\eqref{Submultiplicative}] By \eqref{*Submultiplicative},
\begin{align*}
\|ab\|&=\sqrt{\lceil abb^*a^*\rceil\vee\lceil b^*a^*ab\rceil}\\
&\leq\sqrt{\lceil aa^*\rceil\lceil bb^*\rceil\vee\lceil b^*b\rceil\lceil a^*a\rceil}\\
&\leq\sqrt{(\lceil aa^*\rceil\vee\lceil a^*a\rceil)(\lceil bb^*\rceil\vee\lceil b^*b\rceil)}\\
&=\|a\|\|b\|.
\end{align*}

\item[\eqref{2Subadditive}] Take any $a,b\in A$ and note that
\[(a+b)(a+b)^*=aa^*+ab^*+ba^*+bb^*=2(aa^*+bb^*)-(a-b)(a-b)^*.\]
As $\lceil-p\rceil=0$, for all $p\in|A|^2\subseteq A_\Sigma$, \eqref{Subadditive} yields
\[\lceil(a+b)(a+b)^*\rceil\leq\lceil2(aa^*+bb^*)\rceil\leq2(\lceil aa^*\rceil+\lceil bb^*\rceil)\leq2(\|a\|^2+\|b\|^2).\]
Likewise, $\lceil(a+b)^*(a+b)\rceil\leq2(\|a\|^2+\|b\|^2)$ so
\[\|a+b\|^2\leq2(\|a\|^2+\|b\|^2)\leq2(\|a\|^2+2\|a\|\|b\|+\|b\|^2)=2(\|a\|+\|b\|)^2\]
and hence $\|a+b\|\leq\sqrt{2}(\|a\|+\|b\|)$.
\end{itemize}
\end{proof}

\begin{rmk}
We can use $\|\cdot\|$ to topologise $A$ in the usual way.  Specifically, letting
\[A_\delta=\{a\in A:\|a\|<\delta\},\]
we see that $A_{\delta/3}+A_{\delta/3}\subseteq A_\delta$, by \eqref{2Subadditive}, so the relations $(\equiv^\delta)_{\delta>0}$ (defined from $\|\cdot\|$ as in \eqref{deltaEquiv}) form a uniformity base and hence a topology where $(a+A_\delta)_{\delta>0}$ is a neighbourhood base at each $a\in A$.  So $O\subseteq A$ is open in this topology iff $O$ is a neighbourhood of each $a\in O$, i.e. $a+A_\delta\subseteq O$, for some $\delta>0$.  Note however that, unlike with seminorms, $A_\delta$ itself may not be open.

Alternatively, one could consider the bona fide seminorm
\[|||a|||=\inf\{\tfrac{1}{n}(\|a_1\|+\cdots+\|a_k\|):na=a_1+\cdots+a_k\}.\]
However, this may not be equivalent to $\|\cdot\|$ and may thus induce a coarser topology.  If one prefers $F$-seminorms to quasiseminorms then there is a similar alternative available that is equivalent to $\|\cdot\|$ \textendash\, see \cite[Theorem 1.2]{KaltonPeckRoberts1984}. 

In any case, we will not be considering the topology induced by $\|\cdot\|$ further, instead focusing on algebraic/order properties of $\lceil\cdot\rceil$ and $\|\cdot\|$.
\end{rmk}

To define $\lceil\cdot\rceil$, we could have used `fractions' of $A_\Sigma$.  Specifically, for $B\subseteq A$, let
\[B/\mathbb{N}=\{a\in A:\exists n\in\mathbb{N}\ (na\in B)\}.\]

\begin{prp}\label{Fractions}
Given a *-ring $A$, $a,b\in A$ and $m,n\in\mathbb{N}$,
\[\lceil a\rceil<m/n\qquad\Rightarrow\qquad mbb^*-nbab^*\in A_\Sigma/\mathbb{N}.\]
Conversely, if $mbb^*-nbab^*\in A_\Sigma/\mathbb{N}$, for all $b\in A$, then $\lceil a\rceil\leq m/n$.
\end{prp}

\begin{proof}
If $\lceil a\rceil<m/n$ then, for any $b\in A$, we have $j,k\in\mathbb{N}$ with $j/k<m/n$ and $jbb^*-kbab^*\in A_\Sigma$.  Thus $jn<km$ so
\[k(mbb^*-nbab^*)=(km-jn)bb^*+n(jbb^*-kbab^*)\in A_\Sigma\]
and hence $mbb^*-nbab^*\in A_\Sigma/\mathbb{N}$.

Conversely, if $mbb^*-nbab^*\in A_\Sigma/\mathbb{N}$ then $kmbb^*-knbab^*\in A_\Sigma$, for some $k\in\mathbb{N}$, and $km/(kn)=m/n$.  So if this holds for all $b\in A$ then $\lceil a\rceil\leq m/n$.
\end{proof}

We can sharpen the relationship between $\lceil\cdot\rceil$ and $\|\cdot\|$ if $A$ is unital and self-adjoint/commuting products in $A_\Sigma$ are again in $A_\Sigma$.  If $A$ is also `symmetric' in the sense that adding $1$ to any *-square yields an invertible element (see \cite[\S1 Exercise 7C]{Berberian1972}) then this can be further sharpened to the familiar C*-norm condition.

\begin{thm}\label{C*norm}
If $A$ is a unital *-ring and $(A_\Sigma A_\Sigma)_\mathrm{sa}\subseteq A_\Sigma$ then
\[\|aa^*\|\leq\lceil aa^*\rceil.\]
If, moreover, $1+|A|^2\subseteq A^{-1}(=\{a\in A:a\textrm{ is invertible}\})$ and $\|a\|<\infty$ then
\[\|aa^*\|=\lceil aa^*\rceil=\|a\|^2.\]
\end{thm}

\begin{proof}
For any $\epsilon>0$, we can find $j,k\in\mathbb{N}$ with $j/k\leq\lceil aa^*\rceil+\epsilon$ and $j-kaa^*\in A_\Sigma$.  Then $aa^*(j-kaa^*)=(j-kaa^*)aa^*\in(A_\Sigma A_\Sigma)_\mathrm{sa}\subseteq A_\Sigma$ so
\begin{align*}
j^2-k^2(aa^*)^2&=j^2-kjaa^*+kjaa^*-k^2(aa^*)^2\\
&=j(j-kaa^*)+kaa^*(j-kaa^*)\\
&\in jA_\Sigma+kA_\Sigma\subseteq A_\Sigma.
\end{align*}
Thus $\|aa^*\|=\sqrt{\lceil(aa^*)^2\rceil}\leq\sqrt{j^2/k^2}=j/k\leq\lceil aa^*\rceil+\epsilon$ so $\|aa^*\|\leq\lceil aa^*\rceil$.

Conversely, for any $\epsilon>0$, take $j,k\in\mathbb{N}$ with $\lceil(aa^*)^2\rceil<j^2/k^2\leq\lceil(aa^*)^2\rceil+\epsilon$.  By \autoref{Fractions}, we have some $n\in\mathbb{N}$ with $n(j^2-k^2aa^*)\in A_\Sigma$ and hence $n^2(j^2-k^2aa^*)\in A_\Sigma$.  Replacing $j$ and $k$ with $jn$ and $jk$, we have $j^2-k^2aa^*\in A_\Sigma$.  If $1+|A|^2\subseteq A^{-1}$, then we can let $b=(j+kaa^*)^{-1}$.  Then
\[b=b1^*=b(b(j+kaa^*))^*=b(j+kaa^*)b^*\in A_\Sigma.\]
Thus $A_\Sigma A_\Sigma\ni(j^2-(kaa^*)^2)b=(j-kaa^*)(j+kaa^*)b=j-kaa^*\in A_\mathrm{sa}$
so $j-kaa^*\in A_\Sigma$ and hence $\lceil aa^*\rceil\leq j/k\leq\sqrt{\lceil(aa^*)^2\rceil+\epsilon}$.  As $\epsilon$ was arbitrary, this shows that $\lceil aa^*\rceil\leq\sqrt{\lceil(aa^*)^2\rceil}=\|aa^*\|$.

By \eqref{*Submultiplicative}, we then get
\[\lceil aa^*\rceil^2=\|aa^*\|^2=\lceil aa^*aa^*\rceil\leq\lceil aa^*\rceil\lceil a^*a\rceil,\]
and hence $\lceil aa^*\rceil\leq\lceil a^*a\rceil$.  Likewise $\lceil a^*a\rceil\leq\lceil aa^*\rceil$ so $\|a\|^2=\lceil aa^*\rceil=\lceil a^*a\rceil$.
\end{proof}

\subsection{Order Structure}\label{OrderStructure}

To go further we need to make some further assumptions about our *-ring $A$.  Specifically, we consider $A$ satisfying the following conditions, for all $a\in A$.
\begin{align}
\label{Archimedean}\tag{Archimedean}\lceil-a\rceil=0\qquad&\Rightarrow\qquad\lceil a\rceil<\infty.\\
\label{NoInfinitesimals}\tag{Infinitesimal-Free}\lceil a\rceil=\lceil-a\rceil=0\qquad&\Rightarrow\qquad a=0.\\
\label{Proper}\tag{Proper}aa^*=0\qquad&\Rightarrow\qquad a=0.
\end{align}

\begin{ass}\label{*RingAssumptions}
\textbf{$A$ is an Archimedean infinitesimal-free proper *-ring.}
\end{ass}

Define a relation $\leqslant$ on $A$ by
\[a\leqslant b\qquad\Leftrightarrow\qquad\lceil a-b\rceil=0.\]

\begin{rmk}
In \cite[\S13 Definition 8]{Berberian1972}, the ordering is instead defined by $a-b\in A_\Sigma$.  In $\mathbb{Z}$, for example, both the $A_\Sigma$-order and $\leqslant$ above coincide with the usual ordering, while in $2\mathbb{Z}$, only $\leqslant$ gives the usual ordering.  Incidentally, \eqref{Proper} here is from \cite[\S2 Definition 1]{Berberian1972} (while in commutative ring theory one would instead talk about `reduced' rings).
\end{rmk}

Every C*-algebra is Archimedean, infinitesimal-free and proper, and in this case $\leqslant$ is the usual ordering on $A_\mathrm{sa}$.  In general *-rings, these assumptions suffice to imply that $\leqslant$ has many of the same important properties, as we proceed to show.  We start with the basics.

\begin{prp}
The relation $\leqslant$ is a partial order such that, for all $a,b,c\in A$,
\[a\leqslant b\qquad\Rightarrow\qquad a+c\leqslant b+c\quad\textrm{and}\quad cac^*\leqslant cbc^*.\]
\end{prp}

\begin{proof}
As $\lceil a-a\rceil=\lceil0\rceil$, $\leqslant$ is reflexive.  If $a\leqslant b\leqslant c$ then \eqref{Subadditive} yields
\[\lceil a-c\rceil\leq\lceil a-b\rceil+\lceil b-c\rceil=0.\]
so $a\leqslant c$, i.e. $\leqslant$ is transitive and hence a preorder.  And if $a\leqslant b\leqslant a$ then $\lceil a-b\rceil=\lceil b-a\rceil=0$ so $a=b$, by \eqref{NoInfinitesimals}, i.e. $\leqslant$ is antisymmetric and hence a partial order.

If $a\leqslant b$ then $\lceil(a+c)-(b+c)\rceil=\lceil a-b\rceil=0$, i.e. $a+c\leq b+c$, and $\lceil cac^*-cbc^*\rceil=\lceil c(a-b)c^*\rceil\leq\lceil cc^*\rceil\lceil a-b\rceil=0$, as $\lceil cc^*\rceil<\infty$ by \eqref{Archimedean}, i.e. $cac^*\leqslant cbc^*$.
\end{proof}

We can also describe $\lceil\cdot\rceil$ in terms of $\leqslant$.
\begin{prp}\label{Fractions+}
For any $a,b\in A$ and $m,n\in\mathbb{N}$,
\[\lceil a\rceil\leq m/n\qquad\Rightarrow\qquad nbab^*\leqslant mbb^*.\]
Conversely, if $A$ is unital and $na\leqslant m1$ then $\lceil a\rceil\leq m/n$.
\end{prp}

\begin{proof}
If $\lceil a\rceil\leq m/n$ then, for any $\epsilon>0$ and $c\in A$, we can find $j,k\in\mathbb{N}$ such that $m/n\leq j/k\leq m/n+\epsilon$ and $jcbb^*c^*-kcbab^*c^*\in A_\Sigma$.  By \eqref{Archimedean}, $\lceil bb^*\rceil<g$, for some $g\in\mathbb{N}$, so \autoref{Fractions} yields $h(gcc^*-cbb^*c^*)\in A_\Sigma$, for some $h\in\mathbb{N}$.  Thus
\begin{align*}
&hg(jn-km)cc^*+hkc(mbb^*-nbab^*)c^*\\
=\ &h(jn-km)(gcc^*-cbb^*c^*)+hn(jcbb^*c^*-kcbab^*c^*)\in A_\Sigma
\end{align*}
Note $hg(jn-km)/(hk)=gn(j/k-m/n)=gn\epsilon$.  As $c$ was arbitrary, this shows that $\lceil nbab^*-mbb^*\rceil\leq gn\epsilon$.  As $\epsilon$ was arbitrary, $nbab^*\leqslant mbb^*$.

Conversely, if $A$ is unital and $na\leqslant m1$ then, for any $k\in\mathbb{N}$, $\lceil na-m1\rceil<1/k$ so \autoref{Fractions} yields $j\in\mathbb{N}$ with
\[(j+km)1-kna=j(1-k(na-m1))\in A_\Sigma.\]
As $(j+km)/(kn)=j/(kn)+m/n$ and $k$ was arbitrary, $\lceil a\rceil\leq m/n$.
\end{proof}

Denote the positive elements w.r.t. $\leqslant$ within any $B\subseteq A$ by
\[B_+=\{a\in B:a\geqslant0\}=\{a\in B:\lceil-a\rceil=0\}.\]
Note \eqref{NoInfinitesimals} is saying $A_+$ has no non-trivial additive subgroups, i.e.
\[A_+\cap-A_+=\{0\}.\]
In particular, $A_\Sigma\cap A_\Sigma=\{0\}$, which yields the following.

\begin{prp}\label{A+Asa}
For any $a\in A$ and $n\in\mathbb{N}$,
\[\label{TorsionFree}\tag{Torsion-Free}na=0\qquad\Rightarrow\qquad a=0.\]
Consequently, if $A_+=A_\Sigma/\mathbb{N}$ or if $A$ is unital then $A_+\subseteq A_\mathrm{sa}$.
\end{prp}

\begin{proof}
If $na=0$ then $naa^*=0$ and hence $aa^*=(1-n)aa^*\in A_\Sigma\cap-A_\Sigma=\{0\}$.  By \eqref{Proper}, $a=0$, proving \eqref{TorsionFree}.

In particular, $na=na^*$ implies $a=a^*$, i.e. $A_\mathrm{sa}/\mathbb{N}\subseteq A_\mathrm{sa}$.  Thus if $A_+=A_\Sigma/\mathbb{N}$ then $A_+\subseteq A_\mathrm{sa}/\mathbb{N}$, as $A_\Sigma\subseteq A_\mathrm{sa}$.  Or if $A$ is unital and $a\in A_+$ or even $\lceil-a\rceil<\infty$ then $na\in A_\Sigma-\mathbb{N}1\subseteq A_\mathrm{sa}$, for some $n\in\mathbb{N}$, and hence $a\in A_\mathrm{sa}/\mathbb{N}\subseteq A_\mathrm{sa}$.
\end{proof}

Next we examine properties of $\prec$ on $A$ and its unit ball.

\begin{prp}\label{BallProperties}
Take any $a,b,c\in A$.
\begin{alignat}{7}
\label{a*aprecb}&&&&a^*a&\prec b&&\Rightarrow& a&\prec b.\\
\label{a<pn}&\textrm{If}& b^*=b\geqslant0\quad\textrm{then}&&a&\prec b^n\quad&&\Rightarrow&a&\prec b.\\
\label{prec=>leq}&\textrm{If}& \|a\|\leq1\quad\textrm{then}&&a^*a&\leqslant b^*b&&\Leftarrow\qquad&a&\prec b.\\
\label{aba*=aa*}&\textrm{If}& \|b\|\leq1\quad\textrm{then}&&\quad aba^*&=aa^*\quad&&\Rightarrow&a&\prec b.\\
\label{Ballprec}&\textrm{If}& \|b\|\leq1\quad\textrm{then}&& a&\prec b^*&&\Rightarrow& a&\prec b.\\
\label{BallProduct}&\textrm{If}\quad& \|b\|,\|c\|\leq1\quad\textrm{then}&& a&\prec bc&&\Rightarrow&a&\prec bb^*.\\
\label{precleq}&\textrm{If}& \lceil c\rceil,\|c\|\leq1\quad\textrm{then}&&\quad a\prec b&\leqslant c\quad&&\Rightarrow& a&\prec c.
\end{alignat}
\end{prp}

\begin{proof}\
\begin{itemize}
\item[\eqref{a*aprecb}] If $a^*a\prec b$, i.e. $a^*a=a^*ab$, then
\begin{align*}
(a-ab)^*(a-ab)&=a^*a-a^*ab-b^*a^*a-b^*a^*ab\\
&=a^*a-a^*a-a^*a+a^*a\\
&=0.
\end{align*}
By \eqref{Proper}, $a-ab=0$ so $a=ab$, i.e. $a\prec b$.\\

\item[\eqref{a<pn}] First note that powers of $b$ are positive because $b^{2n}=b^nb^n\in|A|^2\subseteq A_+$ and $b^{2n+1}=b^nbb^n\geqslant b^n0b^n=0$.  If we also have $a=ab^n$ then
\begin{align*}
\qquad0&\leqslant(a-ab)(a-ab)^*\\
&\leqslant(a-ab)(a-ab)^*+(a-ab)b(a-ab)^*+\cdots+(a-ab)b^n(a-ab)^*\\
&=(a-ab^n)(a-ab)^*\\
&=0.
\end{align*}
As $\leqslant$ is antisymmetric, $(a-ab)(a-ab)^*=0$.  By \eqref{Proper}, $a-ab=0$.\\

\item[\eqref{prec=>leq}] If $\|a\|\leq1$ and $a\prec b$ then $a^*a=b^*a^*ab\leqslant b^*b$.\\

\item[\eqref{aba*=aa*}] If $aba^*=aa^*$, taking adjoints yields $ab^*a^*=aa^*$.  If $\|b\|\leq1$ too then
\begin{align*}
0&\leqslant(a-ab)(a-ab)^*\\
&=aa^*-ab^*a^*-aba^*+abb^*a^*\\
&=abb^*a^*-aa^*\\
&\leqslant0.
\end{align*}
So again $a\prec b$, as $\leqslant$ is antisymmetric and \eqref{Proper} holds.\\

\item[\eqref{Ballprec}] If $a\prec b^*$, i.e. $ab^*=a$, then $ab^*ba^*=aa^*$.  If $\|b\|\leq1$ too then $\|bb^*\|\leq\|b\|^2\leq1$ so \eqref{aba*=aa*} yields $a=ab^*b=ab$, as $ab^*=a$, i.e. $a\prec b$.\\

\item[\eqref{BallProduct}] If $\|c\|\leq1$ then $bcc^*b^*\leqslant bb^*$.  If $a\prec bc$ and $\|b\|\leq1$ too then
\[aa^*=abcc^*b^*a^*\leqslant abb^*a^*\leqslant aa^*.\]
As $\leqslant$ is antisymmetric, $abb^*a^*=aa^*$ so $a\prec bb^*$, by \eqref{aba*=aa*}.\\

\item[\eqref{precleq}] If $\lceil c\rceil\leq1$ and $a\prec b\leqslant c$ then $aa^*=aba^*\leqslant aca^*\leqslant aa^*$.  As $\leqslant$ is antisymmetric, $aca^*=aa^*$ so if $\|c\|\leq1$ too then $a\prec c$, by \eqref{aba*=aa*}.
\end{itemize}
\end{proof}

Under further conditions, like in \autoref{C*norm}, we have analogous properties for the \emph{orthogonality} relation $\perp$ defined on $A$ by
\[a\perp b\qquad\Rightarrow\qquad ab^*=0.\]
Note $\perp$ is a symmetric relation, as $ab^*=0$ iff $ba^*=(ab^*)^*=0^*=0$.

\begin{prp}\label{perpResults}
If $p\in A_\Sigma/\mathbb{N}$ or if $p\in A_+\supseteq(A_\Sigma A_\Sigma)_\mathrm{sa}$ and $A$ is unital then
\begin{align}
\label{PositiveProper}apa^*=0\qquad&\Rightarrow\qquad a\perp p.\\
\label{leqprec}p\leqslant a\prec b\qquad&\Rightarrow\qquad p\prec b.\\
\label{leqperp}p\leqslant a\perp b\qquad&\Rightarrow\qquad p\perp b.\\
\label{ProductWedge}p\leqslant a,b\quad\textrm{and}\quad a\perp b\qquad&\Rightarrow\qquad p=0.\\
\intertext{If $p\in A_\Sigma/\mathbb{N}$ and $(A_\Sigma A_\Sigma)_\mathrm{sa}\subseteq A_+$ then}
\label{+perp}ap+pa=0\qquad&\Rightarrow\qquad a\perp p.\\
\label{SquareCommutant}ap^2=p^2a\qquad&\Rightarrow\qquad ap=pa.
\end{align}
\end{prp}

\begin{proof}
Note $p=p^*$, as $p\in A_\Sigma/\mathbb{N}\subseteq A_\mathrm{sa}$ or, if $A$ is unital, $p\in A_+\subseteq A_\mathrm{sa}$.

\begin{itemize}
\item[\eqref{PositiveProper}]
If $p\in A_\Sigma/\mathbb{N}$ and $apa^*=0$ then $0=napa^*=a(b_1b_1^*+\cdots+b_kb_k^*)a^*$, for some $n\in\mathbb{N}$ and $b_1,\cdots,b_k\in A$.  For all $j\leq k$, $ab_jb_j^*a^*\in A_\Sigma\cap-A_\Sigma=\{0\}$ and hence $ab_jb_j^*a^*=0$.  Properness then yields $ab_j=0$ so $ab_jb_j^*=0$ and hence $nap=a(b_1b_1^*+\cdots+b_kb_k^*)=0$.  By \eqref{TorsionFree}, $ap=0$.

Now assume $A$ is unital and $p\in A_+\supseteq(A_\Sigma A_\Sigma)_\mathrm{sa}$.  Then \eqref{Archimedean} yields $m1-np\in A_\Sigma$, for some $m,n\in\mathbb{N}$.  As $p\geqslant0$, for any $\epsilon>0$, we have $j,k\in\mathbb{N}$ with $j/k<\epsilon$ and $j1+kp\in A_\Sigma$.  Then
\[\qquad\qquad jm1+k(mp-np^2)=(j1+kp)(m1-np)+jnp\in(A_\Sigma A_\Sigma)_\mathrm{sa}+A_+\subseteq A_+,\]
i.e. $k(np^2-mp)\leqslant jm1$.  As $jm/k<\epsilon m$ and $\epsilon$ was arbitrary, \autoref{Fractions+} yields $\lceil k(np^2-mp)\rceil=0$, i.e. $np^2\leqslant mp$.  Note $A_+\subseteq A_\mathrm{sa}$, by \autoref{A+Asa}, so if $apa^*=0$ then $0\leqslant nap^2a^*\leqslant mapa^*=0$ and hence $ap=0$, by \eqref{TorsionFree} and \eqref{Proper}.

\item[\eqref{leqprec}] Assume $p\leqslant a\prec b$ and let $c=p-bp-b^*p+bb^*p$ so $ac=ap-ap-ab^*p+ab^*p=0$.  Thus $0\leqslant c^*pc\leqslant c^*ac=0$ so \eqref{PositiveProper} yields $0=pc=(p-pb)(p-pb)^*$ and hence $p-pb=0$, by \eqref{Proper}, i.e. $p\prec b$.

\item[\eqref{leqperp}] If $p\leqslant a\perp b$ then $0\leqslant bpb^*\leqslant bab^*=0$ so $bpb^*=0$.  By \eqref{PositiveProper}, $b\perp p$.

\item[\eqref{ProductWedge}] By \eqref{leqperp}, $p\leqslant a\perp b$ implies $p\perp b$ and then $p\leq b\perp p$ implies $p\perp p$, i.e. $pp^*=0$ so $p=0$, by \eqref{Proper}.

\item[\eqref{+perp}] If $(A_\Sigma A_\Sigma)_\mathrm{sa}\subseteq A_+$ and $p\in A_\Sigma/\mathbb{N}$ then $np\in A_\Sigma\subseteq A_\mathrm{sa}$, for some $n\in\mathbb{N}$.  Then $ap+pa=0$ implies
\[a^*anp=-a^*npa\in(A_\Sigma A_\Sigma)_\mathrm{sa}\cap-A_\Sigma\subseteq A_+\cap-A_+=\{0\}.\]
By \eqref{TorsionFree}, $a^*ap=0$ so $pa^*ap=0=pa^*$, by \eqref{Proper}, i.e. $a\perp p$.

\item[\eqref{SquareCommutant}] If $ap^2=p^2a$ then
\[p(ap-pa)=pap-p^2a=pap-ap^2=(pa-ap)p=-(ap-pa)p\]
and hence $p\perp ap-pa$, by \eqref{+perp}.  Thus $pap=p^2a$ and $pa^*p=p^2a^*$.  Also $aa^*p^2=ap^2a^*=p^2aa^*$ so the same argument applied to $aa^*$ instead of $a$ yields $paa^*p=p^2aa^*$.  Thus
\begin{align*}
(ap-pa)(pa^*-a^*p)&=ap^2a^*-apa^*p-papa^*+paa^*p\\
&=ap^2a^*-ap^2a^*-p^2aa^*+p^2aa^*\\
&=0
\end{align*}
and hence $ap-pa=0$, by \eqref{Proper}, i.e. $ap=pa$.
\end{itemize}
\end{proof}

\begin{rmk}
Note \eqref{SquareCommutant} is saying that $p\in p^{2\prime\prime}$.  It follows that $(A_\Sigma A_\Sigma)_\mathrm{sa}\subseteq A_+$ is equivalent to the commutant part of the positive square-root axiom (PSR) from \cite[\S13 Definition 9]{Berberian1972}, at least under \autoref{*RingAssumptions}.  Specifically, if we have $A_\Sigma\subseteq A_\Sigma^2$, i.e. if every $p\in A_\Sigma$ has a square root $\sqrt{p}\in A_\Sigma$, then
\[\forall p\in A_\Sigma\ (\sqrt{p}\in p'')\qquad\Leftrightarrow\qquad(A_\Sigma A_\Sigma)_\mathrm{sa}\subseteq A_+.\]
Indeed, $\Leftarrow$ is immediate from \eqref{SquareCommutant}.  Conversely, for any $p,q\in A_\Sigma$, $pq\in A_\mathrm{sa}$ means $pq=qp$, i.e. $q\in p'$.  Thus $\sqrt{p}\in p''$ implies $\sqrt{p}q=q\sqrt{p}$ from which we get $pq=\sqrt{p}\sqrt{p}q=\sqrt{p}q\sqrt{p}\in A_\Sigma\subseteq A_+$.  As $p$ and $q$ were arbitrary, $(A_\Sigma A_\Sigma)_\mathrm{sa}\subseteq A_+$.  It would interesting to know if the relevant theory in \cite{Berberian1972} is still valid under $(A_\Sigma A_\Sigma)_\mathrm{sa}\subseteq A_+$, even without square roots.
\end{rmk}

We now return to our investigation of Weyl *-semigroups.  To make use \autoref{BallProperties}, we will consider Weyl *-subsemigroups of the unit ball.

\begin{ass}\label{SsubA1}
\textbf{$(S,E)$ is a Weyl *-semigroup where $S$ lies in the unit ball of an Archimedean infinitesimal-free proper *-ring $A$, i.e.}
\[S\subseteq A^1.\]
\end{ass}

First we obtain a converse of \eqref{TripleProduct}.

\begin{prp}
For any $a,b,c\in S$,
\begin{equation}\label{TripleConverse}
a\precsim bcc^*\qquad\Rightarrow\qquad a\precsim b\quad\textrm{and}\quad a\prec cc^*.
\end{equation}
If $b\sim c$ too then $(bcc^*)^\succsim=b^\succsim\cap c^\succsim$.
\end{prp}

\begin{proof}
If $a\precsim bcc^*$ then $a\prec cc^*b^*bcc^*$ so $a\prec cc^*$, by \eqref{BallProduct}.  Also, as $a\sim bcc^*$, $ab^*=acc^*b^*\in E$.  As $|S|^2$ commutes, $ab^*b=acc^*cc^*b^*b=acc^*b^*bcc^*=a$, proving \eqref{TripleConverse}.  If we also have $c\sim b\succsim a$ then $a\sim c$ too and hence $a\precsim c$.  As $a$ was arbitrary, $(bcc^*)^\succsim\subseteq b^\succsim\cap c^\succsim$.  The reverse inclusion is immediate from \eqref{TripleProduct}.
\end{proof}

With this, we can show that the canonical subbasis $(\mathcal{C}_a)_{a\in S}$ of the coset groupoid $\mathcal{C}=\mathcal{C}(S)$ forms an inverse semigroup of open bisections.  The same then applies to the Weyl groupoid $\mathcal{U}=\mathcal{U}(S)$, as it forms an ideal in $\mathcal{C}$.  Note that we already know each $\mathcal{C}_a$ is a bisection, as $\mathcal{C}_a\mathcal{C}_a^*\subseteq\mathcal{C}_{aa^*}\subseteq\mathcal{C}^0$ and $\mathcal{C}_a^*\mathcal{C}_a\subseteq\mathcal{C}_{aa^*}\subseteq\mathcal{C}^0$, by \autoref{IdempotentCosets}.  The canonical subbasis is also closed under taking inverses, as we always have $\mathcal{C}_a^*=\mathcal{C}_{a^*}$.  It only remains to show closure under products.

\begin{thm}\label{EtaleBasis}
For any $a,b\in S$,
\[\mathcal{C}_{ab}=\mathcal{C}_a\mathcal{C}_b.\]
\end{thm}

\begin{proof}
As $\mathcal{C}_a\mathcal{C}_b\subseteq\mathcal{C}_{ab}$, it suffices to show that $\mathcal{C}_{ab}\subseteq\mathcal{C}_a\mathcal{C}_b$.

Take $U\in\mathcal{C}_{ab}$, so we have $c\in U$ with $c\precsim ab$.  In particular, $c\prec b^*a^*ab$ and hence $c\prec b^*b$, by \eqref{BallProduct}.  Thus $c^*c\in\{b^*b\}^\succ\cap(U^*U)^\precsim$ so $T=(Ub^*)^\precsim$ is a coset, by  \autoref{CaAtlas}.  Also note that $cb^*\precsim abb^*$, by \eqref{TripleProduct}, so $cb^*\precsim a$, by \eqref{TripleConverse}, and hence $a\in T$.  Thus $a^*ab\in(T^*U)^\precsim$ and again \eqref{TripleConverse} yields $b\in(T^*U)^\precsim$.  Thus $U=(T(T^*U)^\precsim)^\precsim\in\mathcal{C}_a\mathcal{C}_b$, as required.
\end{proof}

In the usual situation considered in C*-algebras, $E$ is the unit ball of a Cartan subalgebra $C$ of $A$.  This extra additive structure allows us to prove a few further results on bounds of finitely generated subsets.

\begin{prp}\label{JoinAnalog}
If $A$ is unital, $(A_+A_+)_\mathrm{sa}\subseteq A_+$ and $E$ is the unit ball of the *-subring it generates, for any $a,b,c\in S$ with $a,b\precsim c$, we have $d\in S$ with
\[a^\succsim\cup b^\succsim\subseteq d^\succsim\qquad\textrm{and}\qquad a^\precsim\cap b^\precsim\subseteq d^\precsim.\]
\end{prp}

\begin{proof}
Let $p=a^*a$, $q=b^*b$, $e=p+q-pq$ and $d=ce$.  As $a,b\in S\subseteq A^1$, $0\leqslant p,q\leqslant1$ so $0\leqslant(1-p),(1-q)\leqslant1$ and hence $0\leqslant(1-p)(1-q)\leqslant1$, as $(A_+A_+)_\mathrm{sa}\subseteq A_+$.  As $e=1-(1-p)(1-q)$, it follows that $0\leqslant e\leqslant1$ and hence $e^2\leqslant1$, i.e. $\|e\|\leq1$.  Thus $e\in A^1\cap(E+E-E)\subseteq E$ and hence $d=ce\in S$.

As $a,b\prec c^*c$, it follows that $p,q\prec c^*c$ and hence $e\prec c^*c$ so $d^*d=ec^*ce=e^2$.  If $f\precsim a\precsim c$ then $fd^*=fec^*=fc^*cec^*\in EcEc^*\subseteq E$ and $f\prec a^*a=p$ so $fe=f(p+q-pq)=f+fq-fq=f$ and hence $fd^*d=fe^2=f$, i.e. $f\precsim d$.  Likewise, $f\precsim b$ implies $f\precsim d$ so $a^\succsim\cup b^\succsim\subseteq d^\succsim$.

If $a,b\precsim f$ then $cpf^*=ca^*af^*\in EE\subseteq E$, $cqf^*=cb^*bf^*\in EE\subseteq E$ and $cpqf^*=ca^*ab^*bf^*\in EEE\subseteq E$ so
\[df^*=cef^*=c(p+q-pq)f^*\in A^1\cap(E+E-E)\subseteq E.\]
As $a,b\precsim f$, it also follows that $p,q\prec f^*f$ so $e=p+q-pq\prec f^*f$ and hence $d=ce\prec f^*f$, thus $d\precsim f$.  As $f$ was arbitrary, this shows that $a^\precsim\cap b^\precsim\subseteq d^\precsim$.
\end{proof}

When $C$ is a commutative *-subring, we have a further *-subring $L=C_\mathrm{sa}$ which still contains $|S|^2$.  When $C$ is a C*-algebra, $L$ has other nice properties not shared by $C$, e.g. $L$ is a lattice with respect to $\leqslant$ and consequently generated by its positive elements, i.e. $L=C_+-C_+$.  For the next results, and those on the next section, this is what we really require, namely a nice subset $L$ containing $|S|^2$.  At first we do not even require $L$ to be commutative, although this often follows automatically.

\begin{prp}
If $(A_+A_+)_\mathrm{sa}\subseteq A_+$ and $L\subseteq(L_+-L_+)\cap\mathbb{Z}E_\mathrm{sa}$ then $L\subseteq L'$.
\end{prp}

\begin{proof}
As $L\subseteq L_+-L_+$ and $L\subseteq\mathbb{Z}E$, to show $L\subseteq L'$ it suffices to show $L_+\subseteq E'$.  But for every $p\in L_+\subseteq\mathbb{Z}E_\mathrm{sa}$, we know that $p^2\in|S|^2\subseteq E'$.  This means that, for all $e\in E$, $p^2e=ep^2$ and hence $pe=ep$, by \eqref{SquareCommutant} (or at least the same result and proof with $A_+$ replacing $A_\Sigma$), i.e. $p\in E'$, as required.
\end{proof}

Also, instead of requiring $L\subseteq\mathbb{Z}E_\mathrm{sa}$, it will suffice for $L$ to be `hereditary' in $E$ or even $S$, meaning that any positive element of $L$ below some element of $|S|^2$ must lie within $S$.

\begin{dfn}
An additive subgroup $L$ of $A$ is \emph{$S$-hereditary} if
\begin{equation}\label{LHereditary}
\tag{$S$-Hereditary}|S|^2\subseteq L\qquad\textrm{and}\qquad L_+\cap|S|^{2\geqslant}\subseteq S.
\end{equation}
\end{dfn}

\begin{prp}\label{Trapping}
Take $S$-hereditary $L\subseteq A$ and $p_1,\cdots,p_n,q_1,\cdots,q_n\in|S|^2$ with $p_k\prec q_k$, for all $k\leq n$.  For any $m$ with $1\leq m\leq n$, we have $r\in|S|^2$ with
\[p_1^\succ\cap\cdots\cap p_m^\succ\cap q_{m+1}^\perp\cap\cdots\cap q_n^\perp\subseteq r^\succ\quad\textrm{and}\quad r\in q_1^\succ\cap\cdots\cap q_m^\succ\cap p_{m+1}^\perp\cap\cdots p_n^\perp.\]
\end{prp}

\begin{proof}
Let
\[a=\prod_{1\leq k\leq m}p_k,\quad b=\prod_{m<k\leq n}(a^2-aq_ka)^2\quad\textrm{and}\quad r=b^2.\]
Whenever $m<k\leq n$, $q_k\in|S|^2$ and $S\subseteq A^1$ imply that $aq_ka\leqslant a^2$ and hence $0\leqslant a^2-aq_ka\leqslant a^2$.  As $|S|^2-|S|^2\subseteq L$, this means $a^2-aq_ka\in L_+\cap|S|^{2\geqslant}\subseteq S$.  Thus $b\in S$ and hence $r=b^2\in|S|^2$.

Now assume $c\in p_1^\succ\cap\cdots\cap p_m^\succ\cap q_{m+1}^\perp\cap\cdots\cap q_n^\perp$.  This means $cp_k=c$, for all $k\leq m$, and hence $c=ca=ca^2$.  This also means $q_kc^*=0=cq_k$ and hence
\[c(a^2-aq_ka)=c-cq_ka=c,\]
when $m<k\leq n$, so $cr=cb^2=c$, i.e. $c\prec r$, proving the first inclusion.

As elements of $|S|^2$ commute, $p_k\prec q_k$ implies $a\prec q_k$ and hence $b\prec q_k$, when $1\leq k\leq m$, so $r=b^2\prec q_k$.  Likewise, when $m<k\leq n$, $p_k\prec q_k$ implies \[p_kaq_ka=ap_kq_ka=ap_ka=p_ka^2\]
so $p_k(a^2-aq_ka)=0$ and hence $p_kr=p_kb^2=0$, i.e. $r\in p_k^\perp$.
\end{proof}

We can then extend this result from $|S|^2$ to $S$, as long as our finite collection is bounded above by some element of $S$.

\begin{cor}\label{TrappingCorollary}
Take $S$-hereditary $L\subseteq A$ and $b_1,\cdots,b_n\precsim d$ with $a_k\precsim b_k$, for all $k\leq n$.  For any $m$ with $1\leq m\leq n$, we have $c\in S$ with
\[a_1^\succsim\cap\cdots\cap a_m^\succsim\cap b_{m+1}^\perp\cap\cdots\cap b_n^\perp\subseteq c^\succsim\quad\textrm{and}\quad c\in b_1^\succsim\cap\cdots\cap b_m^\succsim\cap a_{m+1}^\perp\cap\cdots a_n^\perp.\]
\end{cor}

\begin{proof}
By \autoref{Trapping}, we have $r\in|S|^2$ with
\begin{align*}
r\in\ &(b_1^*b_1)^\succ\cap\cdots\cap (b_m^*b_m)^\succ\cap(a_{m+1}^*a_{m+1})^\perp\cap\cdots(a_n^*a_n)^\perp\\
\textrm{and}\quad&(a_1^*a_1)^\succ\cap\cdots\cap(a_m^*a_m)^\succ\cap(b_{m+1}^*b_{m+1})^\perp\cap\cdots\cap(b_n^*b_n)^\perp\subseteq r^\succ.
\end{align*}
Let $c=dr$.  If
\begin{align*}
f&\in a_1^\succsim\cap\cdots\cap a_m^\succsim\cap b_{m+1}^\perp\cap\cdots\cap b_n^\perp\\
&\subseteq(a_1^*a_1)^\succ\cap\cdots\cap(a_m^*a_m)^\succ\cap(b_{m+1}^*b_{m+1})^\perp\cap\cdots\cap(b_n^*b_n)^\perp\\
&\subseteq r^\succ
\end{align*}
then $fc^*c=frd^*dr=fd^*dr=fa_1^*a_1d^*dr=fa_1^*a_1r=fr=f$ and
\[fc^*=frd^*=fa_1^*a_1a_1^*a_1rd^*=fa_1^*a_1ra_1^*a_1d^*\in Ea_1|S|^2a_1^*E\subseteq E|S|^2E\subseteq E\]
so $f\precsim c$.  Also $cb_k^*=drb_k^*=db_k^*b_krb_k^*\in E|S|^2\subseteq E$, for all $k\leq m$, which implies that $c\in b_1^\succsim\cap\cdots\cap b_m^\succsim$, as
$c=dr\in(b_1^*b_1)^\succ\cap\cdots\cap (b_m^*b_m)^\succ$.  Moreover, when $m<k\leq n$, $ca_k^*a_kc^*=ca_k^*a_krd^*=0$ so \eqref{Proper} yields $a_kc^*=0$ and hence $c\in a_{m+1}^\perp\cap\cdots\cap a_n^\perp$ too.
\end{proof}

\subsection{Lattice Structure}\label{LatticeStructure}

As mentioned before \autoref{EtaleBasis}, products in $S$ correspond to products of open bisections in the Weyl groupoid $\mathcal{U}(S)$.  Next, we would like to show that $\precsim$ also has a natural interpretation in $\mathcal{U}(S)$.

Let $\Subset$ denote the `compact containment' relation in a topological space
\[O\Subset N\qquad\Leftrightarrow\qquad\textrm{every open cover of $N$ has a finite subcover of }O.\]
When the topological space is regular one can verify that
\[O\Subset N\qquad\Leftrightarrow\qquad\overline{O}\textrm{ is compact and }\overline{O}\subseteq N.\]
If the space is also locally compact, it follows that $\Subset$ has `interpolation', i.e.
\[O\Subset N\qquad\Rightarrow\qquad\exists\textrm{ open }M\ (O\Subset M\Subset N).\]
What we will aim to do is go in the reverse direction, first showing that $\precsim$ has interpolation.  Using this, we show that $\precsim$ becomes $\Subset$, from which it follows that the Weyl groupoid is locally compact.

Regarding the first step, it is well known that $\prec$ has a certain degree of interpolation on the positive unit ball of a C*-algebra.  Indeed, applying the continuous functional calculus with certain functions $f$ and $g$, one can ensure that
\[\label{abcInterpolation}\tag{$\prec$-Interpolation}a\prec b\prec c\qquad\Rightarrow\qquad a\prec f(b)\prec g(b)\prec c.\]
As $b$ and $c$ are positive and $b=bc=(bc)^*=cb$, they generate a commutative C*-subalgebra $C$.  A moment's thought reveals that all we really need for \eqref{abcInterpolation} is the lattice structure of $C_\mathrm{sa}$, rather than the full continuous functional calculus.  This is more in keeping with our algebraic approach, so accordingly we make the following standing assumption.

\begin{ass}
\textbf{$(S,E)$ is a Weyl *-semigroup, where $S$ lies in the unit ball $A^1$ of an Archimedean infinitesimal-free proper *-ring $A$, and $L$ is a lattice w.r.t. $\leq$ and an $S$-hereditary commutative *-subring of $A$ such that}
\begin{equation}\label{LAssumptions}
L_+=L_\Sigma/\mathbb{N}\quad\textbf{or}\quad A\textbf{ is unital and }(A_+A_+)_\mathrm{sa}\subseteq A_+.
\end{equation}
\end{ass}

For example, when $A$ is a C*-algebra, we would normally take
\[L=C^*(|S|^2)_\mathrm{sa},\]
i.e. the self-adjoint part of the C*-algebra generated by the *-squares from $S$.  As *-squares from $S$ commute, $L=C^*(|S|^2)_\mathrm{sa}$ is indeed a commutative *-subring of $A$, and it is well-known the self-adjoint elements of a commutative C*-algebra form a lattice in their canonical ordering.

In this case, \eqref{LAssumptions} also holds.  Indeed, for C*-algebra $A$, we also always have $A_+=|A|^2=A_\Sigma/\mathbb{N}$ and $(A_+A_+)_\mathrm{sa}\subseteq A_+$.  So if $L=C^*(|S|^2)_\mathrm{sa}$ then indeed $L_+=L_\Sigma/\mathbb{N}$.  We need \eqref{LAssumptions} to use \autoref{perpResults} and also for the following.

\begin{prp}\label{+times-}
For all $a,b\in L$,
\[a\wedge b=0\qquad\Rightarrow\qquad ab=0.\]
\end{prp}

\begin{proof}
If $a,b\in L$ satisfy $a\wedge b=0$ then certainly $a,b\in L_+$.  If $L_+=L_\Sigma/\mathbb{N}$ then we have $m,n\in\mathbb{N}$ with $ma,nb\in L_\Sigma$.  If $ab\neq0$ then we have $p,q\in |L|^2$ with $p\leqslant ma$, $q\leqslant nb$ and $pq\neq0$.  So $q=cc^*$, for some $c\in L$, and \eqref{Archimedean} yields $j,k\in\mathbb{N}$ with $\lceil cc^*\rceil\leq j$ and $\lceil c^*c\rceil\leq k$ and hence $p^2=cc^*cc^*\leqslant kcc^*=kp$ and $\lceil p^2\rceil=\lceil cc^*cc^*\rceil\leq\lceil c^*c\rceil\lceil cc^*\rceil=jk$.  Likewise, $q^2\leqslant hq$ and $\lceil q^2\rceil\leq gh$, for some $g,h\in\mathbb{N}$.  Thus $pq^2p\leqslant ghp^2\leqslant ghkp\leqslant ghkma$ and $qp^2q\leqslant jkq^2\leqslant hjkq\leqslant hjknb$.  As $L$ is commutative, $pq^2p=qp^2q\neq0$, by proper, so $ghkma\wedge hjknb\neq0$.  But $(L,+)$ is a lattice ordered group so, by \cite[Ch XIII \S4 Thoerem 5]{Birkhoff1967}, $a\wedge b=0$ implies $ea\wedge fb$, for all $e,f\in\mathbb{N}$, thus we have a contradiction.

Now instead assume $A$ is unital and $(A_+A_+)_\mathrm{sa}\subseteq A_+$.  As $L$ and hence $L+\mathbb{Z}1$ is commutative, this means $c\leqslant d$ implies $cp\leqslant dp$, for all $c,d\in L+\mathbb{Z}1$ and $p\in L_+$.  By \eqref{Archimedean}, we have $m,n\in\mathbb{N}$ with $a\leqslant m1$ and $b\leqslant n1$, and hence $nab\leqslant mnb$ and $mab\leqslant mna$.  If $ab\neq0$ then \eqref{TorsionFree} and $(A_+A_+)_\mathrm{sa}\subseteq A_+$ again yields $0\lneqq mnab\leqslant mn^2a,m^2nb$, again contradicting $a\wedge b=0=mn^2a\wedge m^2nb$.
\end{proof}

\begin{rmk}
When $(A_+A_+)_\mathrm{sa}\subseteq A_+$ and hence $L_+L_+\subseteq L_+$, $L$ is an \emph{$l$-ring} in the sense of \cite[Ch XVII \S1]{Birkhoff1967}.  In this case, \autoref{+times-} is saying that $L$ is an \emph{almost $f$-ring}, in the sense of \cite[Ch XVII \S6]{Birkhoff1967}.  As $L_+\subseteq A_+\subseteq A_\mathrm{sa}$ and \eqref{Proper} implies $A_\mathrm{sa}$ has no non-zero nilpotents, $L$ is then automatically a true \emph{$f$-ring}, by \cite[Ch XVII \S6 Lemma 1]{Birkhoff1967}.
\end{rmk}

\begin{prp}\label{Joinprec}
For any $p,q,r,s\in|S|^2$,
\[p,q\prec r\prec s\qquad\Rightarrow\qquad p^\succ\cup q^\succ\subseteq(p\vee q)^\succ,\quad(p\vee q)^\perp\subseteq p^\perp\cap q^\perp\quad\textrm{and}\quad p\vee q\prec r.\]
\end{prp}

\begin{proof}
Take $b,c\in S$ with $q=b^*b$, and $r=c^*c$.  If $q\prec r$ then $b\prec c^*c$, by \eqref{a*aprecb}.  As $\|b\|,\|c\|\leq1$, \eqref{prec=>leq} yields $q=b^*b\leqslant c^*cc^*c\leqslant c^*c=r$.  Likewise, $p\prec r$ implies $p\leqslant r$.  If $r\prec s$ too then $p\vee q\leqslant r\prec s$ so \eqref{leqprec} yields $p\vee q\prec s$.

For any $a\prec p\prec r$, note $(r-p)(s-a^*)=r-a^*-p+a^*=r-p$ and hence $0\leqslant r-p\vee q\leqslant r-p\prec s-a^*$ so $r-p\vee q\prec s-a^*$, by \eqref{leqprec}.  This means
\[r-p\vee q=(r-p\vee q)(s-a^*)=r-a^*-p\vee q+(p\vee q)a^*\]
so $a=a(p\vee q)$, i.e. $a\prec p\vee q$ (actually, if $A$ is unital and $(A_+A_+)_\mathrm{sa}\subseteq A_+$ then \autoref{C*norm} yields $\|p\vee q\|\leq\lceil p\vee q\rceil\leq\lceil r\rceil\leq\|c\|\leq1$, in which case $a\prec p\vee q$ is immediate from $a\prec p\leqslant p\vee q$, by \eqref{precleq}).  Likewise, $a\prec q$ implies $a\prec p\vee q$ so $p^\succ\cup q^\succ\subseteq(p\vee q)^\succ$.

As $p,q\leqslant p\vee q$, \eqref{leqperp} immediately yields $(p\vee q)^\perp\subseteq p^\perp\cap q^\perp$.

As above, $r\prec s$ implies $r\leqslant s$ so $0\leqslant s-r$ and hence $p(s-r)=0=q(s-r)$ implies $p\wedge(s-r)=0=q\wedge(s-r)$, by \eqref{ProductWedge}.  Thus $(p\vee q)\wedge(s-r)=0$, as lattice ordered groups are distributive, by \cite[Ch XIII \S4 Theorem 4]{Birkhoff1967}.  Then \autoref{+times-} yields $0=(p\vee q)(s-r)=p\vee q-(p\vee q)r$, i.e. $p\vee q\prec r$.
\end{proof}

This allows us to prove a version of \autoref{JoinAnalog} (which only works for compatible pairs $c,d$ but does not require $E$ to be the unit ball of a *-subring).

\begin{cor}\label{abprecsimcd}
For any $a,b\in S$ and $c,d\in S^\succsim$,
\[a,b\precsim c,d\ \ \textrm{and}\ \ c\sim d\quad\Rightarrow\quad\exists t\in S\ (a^\succsim\cup b^\succsim\subseteq t^\succsim,\ \ t^\perp\subseteq a^\perp\cap b^\perp\ \ \textrm{and}\ \ t\precsim c,d).\]
\end{cor}

\begin{proof}
Let $p=a^*a$ and $q=b^*b$.  By \autoref{Joinprec}, $p\vee q\prec c^*c,d^*d$.  By \eqref{prec=>leq}, $(p\vee q)^2\leqslant(c^*c)^2$ so $(p\vee q)^2\in L_+\cap|S|^{2\geqslant}\subseteq S$ and we may let $t=c(p\vee q)^4\in S$.  Note $t^*t=(p\vee q)^4c^*c(p\vee q)^4=(p\vee q)^8$ so
\[a^\succsim\cup b^\succsim\subseteq p^\succ\cup q^\succ\subseteq(p\vee q)^\succ\subseteq(t^*t)^\succ.\]
If $f\precsim a\precsim c$ then $ft^*=fc^*ct^*=fc^*c(p\vee q)^4c^*\in E|S|^2\subseteq E$ so $f\precsim t$, as $f\in a^\succsim\subseteq(t^*t)^\succ$.  Likewise $f\precsim b$ implies $f\precsim t$ so $a^\succsim\cup b^\succsim\subseteq t^\succsim$.

By \eqref{Proper}, $t^\perp=(t^*t)^\perp=(p\vee q)^\perp\subseteq p^\perp\cap q^\perp=a^\perp\cap b^\perp$.

Moreover, $tc^*=c(p\vee q)^4c^*\in|S|^2\subseteq E$ and $tc^*c=c(p\vee q)^4c^*c=c(p\vee q)^4=t$ so $t\precsim c$.  And finally, $td^*=c(p\vee q)^4d^*=c(p\vee q)^4c^*cd^*\in|S|^2E\subseteq E$, as $c\sim d$, and again $td^*d=c(p\vee q)^4d^*d=c(p\vee q)^4=t$ so $t\precsim d$ as well.
\end{proof}

As usual, we define positive and negative parts of any $a\in L$ by
\[a_+=0\vee a\qquad\textrm{and}\qquad a_-=0\vee-a.\]

\begin{prp}\label{precPositivePart}
For any $a\in A$ and $b\in L$,
\[a\prec b\qquad\Rightarrow\qquad a\prec b_+.\]
\end{prp}

\begin{proof}
As in any lattice ordered group (see \cite[Ch XIII \S4 Theorem 7 (20)]{Birkhoff1967}), $b_+\wedge b_-=0$ and hence $b_+b_-=0$, by \autoref{+times-}.  Also $b=b_+-b_-$ (see \cite[Ch XIII \S3 (15)]{Birkhoff1967}) so $a\prec b$ implies
\[A_+\ni ab_-^2a^*=a(b_--b_+)b_-a^*=-abb_-a^*=-ab_-a^*\in-A_+.\]
Thus $ab_-^2a^*=0=ab_-$, by \eqref{Proper}, and hence $a=ab=ab_+$, i.e. $a\prec b_+$.
\end{proof}

\begin{prp}
For any $p,q\in L_+$, if $p\prec q$ then we have $r,s\in L_+$ with
\begin{align}
\label{precInterpolation2}p^\succ\subseteq r^\succ,\quad p^\perp\subseteq s^\perp,\quad r\prec s\prec q\quad\textrm{and}\quad&s\leqslant q.\\
p\leqslant q\qquad\Rightarrow\qquad &r\leqslant q.\\
p\neq0\qquad\Rightarrow\qquad &r\neq0.
\end{align}
\end{prp}

\begin{proof}
Let $r=(2p-q)_+$ and $s=2p\wedge q$.  If $a\prec p\prec q$, $a(2p-q)=2a-a=a$, i.e. $a\prec 2p-q$ and hence $a\prec(2p-q)_+=r$, by \autoref{precPositivePart}, so $p^\succ\subseteq r^\succ$.

If $p\perp a$ then $s=2p\wedge q\leqslant 2p\perp a$ so $s\perp a$, by \eqref{leqperp}, showing that $p^\perp\subseteq s^\perp$.

As $r,s\leqslant2p\prec q$, \eqref{leqprec} yields $r,s\prec q$ so
\[r-rs=r(q-s)=r(q-2p\wedge q)=r((q-2p)\vee0)=(2p-q)_+(2p-q)_-=0,\]
by \autoref{+times-}, i.e. $r\prec s$.

Also note $p\leqslant q$ implies $2p-q\leqslant 2q-q\leqslant q$ and hence $r=(2p-q)_+\leqslant q$.

As it stands, we could have $r=0$ even when $p\neq0$, but then $r=(2p-q)_+=0$ implies $2p\leqslant q$ and $p^\succ\subseteq r^\succ=0^\succ=\{0\}$.  However, $p\neq0$ implies $np\nleqslant q$, for some $n\in\mathbb{N}$ \textendash\, otherwise $np\leqslant q$ would imply $n\lceil p\rceil\leq\lceil q\rceil$, for all $n\in\mathbb{N}$, and hence $\lceil p\rceil=0=\lceil-p\rceil$, contradicting \eqref{NoInfinitesimals}.  Thus we must have $n\in\mathbb{N}$ with $np\leqslant q$ but $2np\nleqslant q$.  Also $(np)^\perp=p^\perp$, by \eqref{TorsionFree}, and $p^\succ=\{0\}\subseteq\{np\}^\succ$ so we can replace $p$ with $np$ at the start to ensure that $r\neq0$.
\end{proof}

This interpolation property can then be extended from $L_+$ to $S$.

\begin{cor}\label{SInterpolation}
For any $a,b\in S$ with $a\precsim b$, we have $c,d\in S$ with
\begin{align*}
a^\succsim\subseteq c^\succsim,\quad a^\perp\subseteq d^\perp\quad\textrm{and}\quad&c\precsim d\precsim b.\\
a\neq0\quad\ \Rightarrow\ \quad&c\neq0.
\end{align*}
\end{cor}

\begin{proof}
Let $p=a^*a$ and $q=b^*b$ and take $r,s\in L_+(\subseteq A_\mathrm{sa})$ as in \eqref{precInterpolation2}.  Note $a\prec b^*b$, as $a\precsim b$, and hence $p=a^*a\leqslant b^*bb^*b\leqslant b^*b=q$, by \eqref{prec=>leq}.  Thus $r,s\leqslant q$ and hence $r,s\in L_+\cap|S|^{2\geqslant}\subseteq S$ so $r^2,s^2\in |S|^2\subseteq E$.  Let
\[c=br^2\in S\qquad\textrm{and}\qquad d=bs^2\in S.\]
Then $c^*c=r^2b^*br^2=r^2qr^2=r^4$ and $d^*d=s^2b^*bs^2=s^2qs^2=s^4$ so
\[(a^*a)^\succ\subseteq(c^*c)^\succ,\quad c\prec d^*d\quad\textrm{and}\quad d\prec b^*b.\]
If $f\precsim a\precsim b$ then $fc^*=fr^2b^*=fb^*br^2b^*\in E|S|^2\subseteq E$ and $f\in(a^*a)^\succ\subseteq(c^*c)^\succ$ so $f\precsim c$.  As $f$ was arbitrary, $a^\succsim\subseteq c^\succsim$.  As $p\in Sa$ and $d\in Ss$, $a^\perp\subseteq p^\perp\subseteq s^\perp\subseteq d^\perp$.  Also $cd^*=br^2s^2b^*=br^2b^*\in|S|^2\subseteq E$ so $c\precsim d$, as we already know $c\prec d^*d$.  Likewise, $db^*=bs^2b^*\in|S|^2\subseteq E$ so $d\precsim b$, as we already know $d\prec b^*b$.  Lastly, by \eqref{Proper}, if $a\neq0$ then $p=a^*a\neq0$ so $r\neq0$ and hence $c^*c=r^4\neq0\neq c$.
\end{proof}

\begin{prp}\label{Uperp}
For any $U\in\mathcal{U}(S)$,
\[(U^\succsim\setminus U)^{\succsim\succsim}\subseteq U^\perp\subseteq S\setminus U.\]
\end{prp}

\begin{proof}
Take $a\in(U^\succsim\setminus U)^{\succsim\succsim}$ so $a\precsim b\precsim c\precsim u\in U$, for some $c\notin U$.  By \autoref{SInterpolation}, we have a sequence $(b_n)\subseteq S$ with $b^\succsim\subseteq b_n^\succsim$ and $b_{n+1}\precsim b_n\precsim c$, for all $n\in\mathbb{N}$.  Let
\[T=\{t\in S:\exists v\in U\ \exists n\in\mathbb{N}\ (v^\succsim\cap b_n^\succsim\subseteq t^\succsim)\}.\]
Certainly $T^\precsim\subseteq T$.  Also, for any $s,t\in T$, we have $v,w\in U$ and $j,k\in\mathbb{N}$ with $v^\succsim\cap b_j^\succsim\subseteq s^\succsim$ and $w^\succsim\cap b_k^\succsim\subseteq t^\succsim$.  W.l.o.g. assume $j\leq k$.  As $U$ is a filter, we have $x,y\in U$ with $x\precsim y\precsim u,v,w$.  As $b_{k+1}\precsim b_k\precsim c\precsim u$, \autoref{TrappingCorollary} yields $z\precsim b_k,y$ with $x^\succsim\cap b_{k+1}^\succsim\subseteq z^\succsim$ and hence $z\precsim s,t$ and $z\in T$.  Thus $T$ is directed and hence a filter.  We also immediately see that $U\subsetneqq T$, as $c\in T\setminus U$, so $T=S$, as $U$ is an ultrafilter.  In particular, $0\in T$ so we have $v\in U$ and $n\in\mathbb{N}$ with $v^\succsim\cap b_n^\succsim\subseteq0^\succsim$.  Take $w\in U$ with $w\precsim u,v$.  As $a\precsim u\sim w$, we have $a\sim w$ so $aw^*w\precsim b,v$, by \eqref{EBelow}, as $a\precsim b$ and $w\precsim v$.  As $b^\succsim\subseteq b_n^\succsim$, it follows that $aw^*w\in v^\succsim\cap b_n^\succsim\subseteq0^\succsim$ so $aw^*w=aw^*w00=0$.  Thus $aw^*wa^*=0=wa^*$, by \eqref{Proper}, so $a\in U^\perp$.  As $a$ was arbitrary, this proves $(U^\succsim\setminus U)^{\succsim\succsim}\subseteq U^\perp$.  As filters are cosets, if $v\perp u$, for some $u,v\in U$, then $0=uv^*v\in U$, a contradiction, and hence $U^\perp\subseteq S\setminus U$ too.
\end{proof}

\begin{prp}\label{UltrafilterCovers}
For any $T\subseteq S$ and $a\precsim b$,
\[\mathcal{U}_b\subseteq\bigcup_{t\in T}\mathcal{U}_t\quad\Rightarrow\quad\exists\textrm{ finite }F\subseteq(T^\succsim\cap b^\succsim)^{\succsim\succsim}\ (a^\succsim\cap\bigcap_{f\in F}f^\perp=\{0\})\quad\Rightarrow\quad\mathcal{U}_a\subseteq\bigcup_{t\in T}\mathcal{U}_t.\]
\end{prp}

\begin{proof}
Assume $a^\succsim\cap\bigcap_{f\in F}f^\perp\neq\{0\}$, for all finite $F\subseteq (T^\succsim\cap b^\succsim)^{\succsim\succsim}$.  We aim to construct $U\in\mathcal{U}_b\setminus\bigcup_{t\in T}\mathcal{U}_t$.  As $a\precsim b$, repeated applications of \autoref{SInterpolation} yield a sequence $(a_n)\subseteq S$ such that $a^\succsim\subseteq a_n^\succsim$ and $a_{n+1}\precsim a_n\precsim b$, for all $n\in\mathbb{N}$.  Let
\[U=\{u\in S:\exists n\in\mathbb{N}\ \exists\textrm{ finite }F\subseteq(T^\succsim\cap b^\succsim)^{\succsim\succsim}\ (a_n^\succsim\cap\bigcap_{f\in F}f^\perp\subseteq u^\succsim)\}.\]
Note $\{0\}\neq a^\succsim\cap\bigcap_{f\in F}f^\perp\subseteq a_n^\succsim\cap\bigcap_{f\in F}f^\perp\subseteq u^\succsim$ implies $u\neq0$ so $0\notin U$.  Also note $U^\precsim\subseteq U$.  Moreover, for any $u,v\in U$, we have finite $F,G\subseteq(T^\succsim\cap b^\succsim)^{\succsim\succsim}$ and $j,k\in\mathbb{N}$ with $a_j^\succsim\cap\bigcap_{f\in F}f^\perp\subseteq u^\succsim$ and $a_k^\succsim\cap\bigcap_{g\in G}g^\perp\subseteq v^\succsim$.  W.l.o.g. assume $j\leq k$.  By \autoref{SInterpolation}, we have finite $H\subseteq(F\cup G)^\precsim\cap (T^\succsim\cap b^\succsim)^{\succsim\succsim}$ with $F\cup G\subseteq H^\succsim$.  By \autoref{TrappingCorollary}, we have $w\in a_k^\succsim\cap\bigcap_{e\in F\cup G}e^\perp$ with $a_{k+1}^\succsim\cap\bigcap_{h\in H}h^\perp\subseteq w^\succsim$ so $w\precsim u,v$ and $w\in U$.  Thus $U$ is directed and hence a filter.  By the Kuratowski-Zorn lemma, $U$ extends to an ultrafilter, which we again denote by $U$.  Note $b\in U$, as $b\succsim a_1$.  If we had $U\cap T\neq\emptyset$ then we would have $x,y\in U\cap(T^\succsim\cap b^\succsim)^{\succsim\succsim}$ with $x\precsim y$.  Then \autoref{TrappingCorollary} yields $z\in a_1^\succsim\cap x^\perp$ with $a_2^\succsim\cap y^\perp\subseteq z^\succsim$ and hence $z\in U$, even though $z\perp x\in U$, a contradiction.  Thus $U\cap T=\emptyset$, i.e. $U\in\mathcal{U}_b\setminus\bigcup_{t\in T}\mathcal{U}_t$, proving the first implication.

For the second, assume we had
\[U\in\mathcal{U}_a\setminus\bigcup_{t\in T}\mathcal{U}_t=\mathcal{U}_a\setminus\bigcup_{t\in T^\succsim}\mathcal{U}_t\subseteq\mathcal{U}_b\setminus\bigcup_{t\in T^\succsim}\mathcal{U}_t,\]
as $a\precsim b$.  It follows that $T^\succsim\cap b^\succsim\subseteq b^\succsim\setminus U\subseteq U^\succsim\setminus U$ and hence \autoref{Uperp} yields $(T^\succsim\cap b^\succsim)^{\succsim\succsim}\subseteq(U^\succsim\setminus U)^{\succsim\succsim}\subseteq U^\perp$.  Thus, for any finite $F\subseteq(T^\succsim\cap b^\succsim)^{\succsim\succsim}\subseteq U^\perp$, we have $c\in U$ with $c\precsim a$ and $c\perp f$, for all $f\in F$, as $U$ is a filter (noting that $c\precsim u\perp f$ implies $cf^*=cu^*uf^*=0$, i.e. $c\perp f$).  As $c\neq0$, this shows that $a^\succsim\cap\bigcap_{f\in F}f^\perp\neq\{0\}$, as required.
\end{proof}

\begin{prp}\label{UaInclusionUb}
For any $a,b\in S$,
\[\forall c\precsim a\ \exists d\in(a^\succsim\cap b^\succsim)^\succsim\ (c^\succsim\cap d^\perp=\{0\})\quad\Leftrightarrow\quad a^{\succsim\succsim}\subseteq b^{\succsim\succsim}\quad\Rightarrow\quad\mathcal{U}_a\subseteq\mathcal{U}_b.\]
The last implication can be reversed if $E$ is the unit ball of the *-subring it generates or if $a,b\in S^\succsim$ and $a\sim b$.
\end{prp}

\begin{proof}
For the first $\Rightarrow$, let $f\in a^{\succsim\succsim}\setminus b^{\succsim\succsim}$.   Take $c\in S$ with $f\precsim c\precsim a$.  If $d\precsim e\precsim a,b$ then $e\sim f$, as $e,f\precsim a$, and hence $f\not\prec e^*e$ (otherwise $f\precsim e\precsim b$, contradicting $f\notin b^{\succsim\succsim}$).  Thus $0\neq f-fe^*e=f(a^*a-e^*e)$ and $(a^*a-e^*e)\in S$, as $e^*e\leqslant a^*aa^*a\leqslant a^*a$ and hence $a^*a-e^*e\in L_+\cap|S|^{2\geqslant}$.  Let $g=f(a^*a-e^*e)^2$ and note $g\neq0$, as $gf^*=f(a^*a-e^*e)^2f^*=(f-fe^*e)(e^*ef^*-f^*\neq0$, by \eqref{Proper}.  Then $gc^*=f(a^*a-e^*e)^2c^*=fc^*c(a^*a-e^*e)^2c^*\in E|S|^2\subseteq E$ and also $gc^*c=f(a^*a-e^*e)^2c^*c=fc^*c(a^*a-e^*e)^2=f(a^*a-e^*e)^2=g$ so $g\precsim c$.  Also $dg^*=d(a^*a-e^*e)^2f^*=0$, as $d\precsim e\precsim a$, so $g\in c^\succsim\cap d^\perp$.  Thus we have shown that $a^{\succsim\succsim}\nsubseteq b^{\succsim\succsim}$ implies $\exists c\precsim a\ \forall d\in(a^\succsim\cap b^\succsim)^\succsim\ (c^\succsim\cap d^\perp\neq\{0\})$, as required.

Conversely, assume $a^{\succsim\succsim}\subseteq b^{\succsim\succsim}$ and take $c\precsim a$.  By two applications of \autoref{SInterpolation}, we have $d\in a^{\succsim\succsim\succsim}=(a^{\succsim\succsim}\cap b^{\succsim\succsim})^\succsim\subseteq(a^\succsim\cap b^\succsim)^\succsim$ such that $c^\succsim\subseteq d^\succsim$ and hence $c^\succsim\cap d^\perp\subseteq d^\succsim\cap d^\perp=\{0\}$.

For the second $\Rightarrow$, note that if $a^{\succsim\succsim}\subseteq b^{\succsim\succsim}$ then, for any $U\in\mathcal{U}_a$, we know that $\emptyset\neq U\cap a^{\succsim\succsim}\subseteq U\cap b^{\succsim\succsim}$ and hence $b\in U$, i.e. $U\in\mathcal{U}_b$.

Lastly, assume $E$ is the unit ball of the *-subring it generates or $a,b\in S^\succsim$ and $a\sim b$.  If $\mathcal{U}_a\subseteq\mathcal{U}_b$ then, for any $c\precsim a$, \autoref{UltrafilterCovers} yields finite $F\subseteq(a^\succsim\cap b^\succsim)^\succsim$ with $c^\succsim\cap\bigcap_{f\in F}f^\perp=\{0\}$.  This means we have finite $G\subseteq a^\succsim\cap b^\succsim$ with $F\subseteq G^\succsim$.  By \autoref{JoinAnalog} or \autoref{abprecsimcd}, we have $h\precsim a,b$ with $F\subseteq G^\succsim\subseteq h^\succsim$.  By \autoref{abprecsimcd} again, we have $d\precsim h\precsim a,b$ with $d^\perp\subseteq\bigcap_{f\in F}f^\perp$ and hence $c^\succsim\cap d^\perp=\{0\}$, as required.
\end{proof}

\begin{thm}\label{UaSubsetUb}
For any $a,b\in S$,
\[\exists c\in S\ (a^{\succsim\succsim}\subseteq c^{\succsim\succsim}\textrm{ and }c\precsim b)\qquad\Rightarrow\qquad\mathcal{U}_a\Subset\mathcal{U}_b.\]
If $E$ is the unit ball of the *-subring it generates or if $a,b\in S^\succsim$ and $a\sim b$ then the converse also holds.
\end{thm}

\begin{proof}
Assume $a^{\succsim\succsim}\subseteq c^{\succsim\succsim}$ and $c\precsim b$.  Given a (basic) open cover $(\mathcal{U}_t)_{t\in T}$ of $\mathcal{U}_b$, \autoref{UltrafilterCovers} yields finite $F\subseteq(T^\succsim\cap b^\succsim)^{\succsim\succsim}$ with
\[a^{\succsim\succsim}\cap\bigcap_{f\in F}f^\perp\subseteq c^{\succsim\succsim}\cap\bigcap_{f\in F}f^\perp\subseteq c^\succsim\cap\bigcap_{f\in F}f^\perp=\{0\}.\]
This implies $a^\succsim\cap\bigcap_{f\in F}f^\perp=\{0\}$ \textendash\, if we had non-zero $g\in a^\succsim\cap\bigcap_{f\in F}f^\perp$ then \autoref{SInterpolation} would yield non-zero $d,e$ with $d\precsim e\precsim a$ and $F\subseteq g^\perp\subseteq e^\perp\subseteq d^\perp$ so $d\in a^{\succsim\succsim}\cap\bigcap_{f\in F}f^\perp\neq\{0\}$.  Taking finite $G\subseteq T$ such that $F\subseteq(G^\succsim\cap b^\succsim)^{\succsim\succsim}$, \autoref{UltrafilterCovers} again yields $\mathcal{U}_a\subseteq\bigcup_{g\in G}\mathcal{U}_g$.  Thus every open cover of $\mathcal{U}_b$ has a finite subcover of $\mathcal{U}_a$, i.e. $\mathcal{U}_a\Subset\mathcal{U}_b$.

Conversely, assume $a,b\in S^\succsim$, $a\sim b$ and $\mathcal{U}_a\Subset\mathcal{U}_b$.  As $\mathcal{U}_b\subseteq\bigcup_{f\in b^\succsim}\mathcal{U}_f$, we have finite $F\subseteq b^\succsim$ with $\mathcal{U}_a\subseteq\bigcup_{f\in F}\mathcal{U}_f$.  By \autoref{abprecsimcd}, we have $c\precsim b$ with $F^\succsim\subseteq c^\succsim$ and hence $\mathcal{U}_a\subseteq\bigcup_{f\in F}\mathcal{U}_f\subseteq\mathcal{U}_c$ so $a^{\succsim\succsim}\subseteq c^{\succsim\succsim}$, by \autoref{UaInclusionUb}.
\end{proof}

We can now finally reach the desired local compactness result.

\begin{cor}\label{WeylLocallyCompact}
The Weyl groupoid $\mathcal{U}(S)$ is locally compact.
\end{cor}

\begin{proof}
For any (basic) open neighbourhood $\mathcal{U}_a$ of $U\in\mathcal{U}(S)$, we have $b\in U$ with $b\precsim a$ and hence $U\in\mathcal{U}_b\Subset\mathcal{U}_a$, by \autoref{UaSubsetUb}.
\end{proof}

\newpage


\bibliographystyle{alphaurl}
\bibliography{maths}

\end{document}